\documentclass[nonblindrev]{write_paper} 
\usepackage{amsmath}
\usepackage{url}
\usepackage{algorithm}
\usepackage{booktabs}
\usepackage{multirow}
\usepackage{algpseudocode}
\usepackage{subfigure}
\usepackage{epsfig}

\makeatletter
\renewcommand\section{\@startsection {section}{1}{\z@}%
                                   {-1ex \@plus -1ex \@minus -.1ex}%
                                   {1 ex \@plus.1ex}%
                                   {\normalfont\large\bfseries}}
\renewcommand\subsection{\@startsection{subsection}{2}{\z@}%
                                     {-1ex\@plus -1ex \@minus -.1ex}%
                                     {1ex \@plus .1ex}%
                                     {\normalfont \normalsize \bfseries}}
\renewcommand\subsubsection{\@startsection{subsubsection}{3}{\z@}%
                                     {-1ex\@plus -1ex \@minus -.1ex}%
                                     {1ex \@plus .1ex}%
                                     {\normalfont\normalsize\bfseries}}
\makeatother

 \renewcommand{\theARTICLETOP}{}
\usepackage{fancyhdr}
\usepackage[T1]{fontenc} \usepackage{palatino}

\OneAndAHalfSpacedXI 
\usepackage{color}
\usepackage{soul}

\DeclareMathOperator{\E}{\mathbb{E}}
\DeclareMathOperator{\R}{\mathbb{R}}
\DeclareMathOperator{\B}{\mathbb{B}}

\usepackage{natbib}
 \bibpunct[, ]{(}{)}{,}{a}{}{,}%
 \def\bibfont{\small}%
\TheoremsNumberedThrough

\EquationsNumberedThrough
\MANUSCRIPTNO{} 

\newtheorem{prop}{\textbf{Proposition}}

\newcommand{\JJ}[1]{{\color{black}#1}}
\newcommand{\eat}[1]{}
\usepackage{graphicx}

\usepackage{endnotes}
\let\footnote=\endnote

\usepackage[symbol]{footmisc}

\newcommand{\bs}[1]{\boldsymbol{#1}}
\newcommand{\ml}[1]{\mathcal{#1}}
\newcommand{\mb}[1]{\mathbb{#1}}

\begin{document}




\TITLE{Approximate Resolution of Stochastic Choice-based Discrete Planning}


	\ARTICLEAUTHORS{%
                    \AUTHOR{Jiajie Zhang}
                    \AFF{Department of Industrial Systems Engineering and Management, National University of Singapore, Republic of Singapore, \EMAIL{jiajiez@u.nus.edu}}

		\AUTHOR{Yun Hui Lin}
		\AFF{Institute of High Performance Computing (IHPC), Agency for Science, Technology and Research (A*STAR), 1 Fusionopolis Way, \#16-16 Connexis, Singapore 138632, Republic of Singapore, \EMAIL{liny@ihpc.a-star.edu.sg} \URL{}}

		\AUTHOR{Gerardo Berbeglia}
		\AFF{Melbourne Business School, The University of Melbourne, \EMAIL{g.berbeglia@mbs.edu} \URL{}}

	} 

\ABSTRACT{%

\noindent \textbf{Abstract.} Stochastic choice-based discrete planning is a broad class of decision-making problems characterized by a sequential decision-making process involving a planner and a group of customers. The firm or planner first decides a subset of options to offer to the customers, who, in turn, make selections based on their utilities of those options. This problem has extensive applications in many areas, including assortment planning, product line design, and facility location. A key feature of these problems is that the firm cannot fully observe the customers' utilities or preferences, which results in intrinsic and idiosyncratic uncertainties. Most works in the literature have studied a specific type of uncertainty, resulting in customized decision models that are subsequently tackled using ad-hoc algorithms designed to exploit the specific model structure. 

In this paper we propose a modeling framework capable of solving this family of sequential problems that works for a large variety of uncertainties. We then leverage an approximation scheme and develop an adaptable mixed-integer linear programming method. To speed up the solution process, we further develop an efficient decomposition approach. We show that our solution framework can yield solutions proven to be (near-)optimal for a broad class of problems. We illustrate this by applying our approach to three classical application problems: constrained assortment optimization and two facility location problems. Through extensive computational experiments, we demonstrate the performance of our approach in terms of both solution quality and computational speed, and provide computational insights. In particular, when we use our method to solve the constrained assortment optimization problem under the Exponomial choice model, it improves the state-of-the-art.

}



\maketitle

%


\vspace{-5mm}

\section{Introduction}

The discrete planning problem (DPP) represents a broad class of decision-making problems in which the planner must select a subset from a finite set of options (e.g., products, services, facilities, schedules). It encompasses a wide range of applications, including vehicle routing, machine scheduling, facility location, and assortment planning, among others. In this paper, we focus on a particular branch of DPP, referred to as \textit{choice-based DPP}. This variant investigates a sequential decision-making process involving two primary stakeholders: a planner and a group of customers. The planner determines a subset of options to offer the customers, who in turn make selections based on specified choice rules. These selections play a crucial role in evaluating the value of the planner's objective function.  Choice-based DPP has extensive applications in areas such as assortment optimization~\citep{talluri2004revenue, abeliuk2016assortment,alfandari2021exact,altekin2021linear}, product line design~\citep{bertsimas2019exact}, facility location planning~\citep{fernandez2017new,haase2014comparison,lamontagne2023optimising,lin2022locating,ljubic2018outer}, and the joint optimization of facility location and service pricing~\citep{kochetov2015comparison,lin2023facility,zambrano2019retail}.

\eat{The core of choice-based DPP revolves around the understanding and anticipation of customer choices. In the literature, the utility-based approach is typically adopted for modeling choices, where utility quantifies an individual's satisfaction or desirability associated with a specific option. Notably, a substantial body of research assumes that  the utility is \textit{deterministic and known} to the planner and that customers exhibit \textit{rational} decision-making behaviors. As a result, customers exclusively select the option with the highest utility, leading to the deterministic first-choice outcome. This assumption enables the formulation of DPP as various mixed-integer linear programs (MILPs), as demonstrated in previous works \citep{mcbride1988integer, belloni2008optimizing, bertsimas2019exact}. Apart from MILP formulations, \citet{lin2024revisiting} formulate a bilevel model and propose an exact branch-and-cut algorithm to solve it. Additionally, approximate algorithms \citep{aouad2018approximability} and dynamic-programming-based approaches \citep{aouad2021assortment} have been proposed to tackle first-choice DPP.}

The core of choice-based DPP revolves around understanding and anticipating customer choices. In the literature, the utility-based approach is typically adopted for modeling choices, where customers are assumed to be \textit{economic} or \textit{rational}, preferring the option with the highest utility. Often, the utility that a customer assigns to an option can be divided into two parts: a component that is intrinsic to the option and common to all customers, and an \emph{idiosyncratic} component that reflects the unique customer preferences. To illustrate, consider the problem of choosing a route between two locations: the intrinsic (also called external) utility may be determined by factors such as travel time and distance for each route, while the idiosyncratic (also known as internal) utility could be linked to the driver's personal knowledge and familiarity with those routes. Since it is not always possible to have direct access to the idiosyncratic utility of its customers, and often, there is also not enough information about the intrinsic utility of the alternatives, firms should make planning decisions under uncertainty. Consequently, the concept of stochastic choice-based DPP becomes highly relevant in practice.

A general framework for modeling the uncertainty about customer preferences involves using \emph{random utility maximization models} (RUMs)~\citep{mcfadden1972conditional,mcfadden1981econometric}. Under a RUM, a firm only knows the probability distribution of the utilities of each alternative. When a customer is faced with a subset of alternatives, the utilities of these alternatives are realized, and she selects the one with the highest utility. Thus, the choice behaviors of customers take on a probabilistic nature. Depending on different assumptions on the probability distribution of the utilities, RUMs have given rise to different models, among which the multinomial logit model \citep{zambrano2019retail, abeliuk2016assortment, davis2013assortment}, the nested logit model \citep{alfandari2021exact, chen2020assortment, gallego2014constrained}, the exponomial choice model \citep{alptekinouglu2016exponomial,aouad2023exponomial} and the mixed logit model~ \citep{haase2014comparison, sen2018conic, shen2009latent} are popular because of their interpretability and choice structure. Additionally, RUMs are able to capture choice models based on different primitives such as the Markov chain choice model and the random consideration set choice model \citep{blanchet2016markov,berbeglia2016discrete,manzini2014stochastic}. Recently, there has been new models that subsume RUMs such as the generalized stochastic preference choice model \citep{berbeglia2018generalized} and the decision forest choice model \citep{chen2022decision, chen2019use}.


Undoubtedly, the uncertainties in choice-based DPP add substantial challenges to solve the well-known assortment optimization problem in revenue management. The prevalent literature on choice-based DPP typically presumes a specific structure of uncertainties and develops ad-hoc algorithms. On one hand, this convention necessitates extensive domain knowledge. On the other hand, the resulting customized models often exhibit limited adaptability. For example, among RUM models, the most widely used is arguably the multinomial logit model. It posits that idiosyncratic uncertainty in the utility follows the Gumbel distribution, leading to a linear fractional choice structure that largely facilitates the design of efficient exact solution approaches. However, when the probability distribution of utility uncertainty shifts to the Normal distribution, the model transitions to a probit choice model, for which there is currently no efficient methodology to solve choice-based DPP. Therefore, it is crucial to develop a more unified and adaptive paradigm capable of modeling and solving choice-based DPP that can accommodate different distributions of utility uncertainty, whether intrinsic or idiosyncratic.

The Sample Average Approximation (SAA) scheme has been used to address the challenges posed by uncertainties. The idea is to sample a finite set of scenarios to approximate the uncertainties according to the associated probability distribution. In each scenario, the utility becomes deterministic. By approximating choice probabilities using the sample mean, several mixed-integer programming (MIP) formulations have been developed. For example, \citet{haase2016maximum} proposed two MIP formulations for a maximum capture facility location problem under RUMs. The second formulation has fewer constraints but requires more samples to achieve a solution quality similar to that of the first one. On this basis, \citet{legault2024model} proposed a more compact maximum covering formulation as well as a partial Benders decomposition approach leveraging the submodularity of the objective. \citet{lamontagne2023optimising} proposed a new greedy heuristic to solve the resulting maximum covering formulation and applied it to a charging station location problem.
\citet{paneque2021integrating} considered a profit maximizing planning problem and proposed a more general mixed-integer linear programming (MILP) formulation leveraging  ``big-M'' constraints related to utilities.
 
 The closest related papers are \citet{legault2024model} where the authors developed a sample-based approximation methodology which could be extended to cover both intrinsic and idiosyncratic uncertainties.  However, their work exclusively focuses on the maximum capture facility location problem, where the objective is to maximize the market share. In this case, the ''reward'' that the firm receives for serving customers at any facility is the same. Based on such a homogeneous reward setting, their SAA model can be visualized as a type of set covering problems. It shows submodularity and allows for scenario reduction, which facilitate the development of a partial Benders decomposition algorithm. We note that their modeling and solution approach is not applicable to applications with heterogeneous rewards, such as assortment optimization where products could have different selling prices. Moreover, it cannot address problems that does not employ RUMs to model the choices of customers, such as the partially binary choice rule~\citep{fernandez2017new}. \citet{paneque2018lagrangian, paneque2022lagrangian} studied a more general choice-based optimization problem and proposed a Lagrangian relaxation method when solving the approximate SAA model. While this approach has broader applicability, it cannot guarantee an optimal solution to the approximate model because of the duality gap. This gap could potentially magnify the approximation error within the SAA method and introduce additional issues when attempting to validate the solution quality through statistical inferences. Furthermore, the above works only consider a deterministic reward, which is independent of uncertainties. However, as we will see in later sections, the reward could be closely related to uncertainties and even correlated with utilities in some applications.

Motivated by the above challenges and limitations, this paper develops a unified and adaptable approach for streamlining the modeling and algorithmic resolution of stochastic choice-based DPP. The contributions can be summarized as follows:
\begin{itemize}
\item Unified model. We present a unified model for the stochastic choice-based DPP, which is not limited to specific problem settings and is applicable to both idiosyncratic (or internal) and intrinsic (or external) uncertainties. To illustrate its general applicability, we provide three examples of applications, i.e., the classical cardinality constrained assortment optimization problem (CAOP) under random utilities, and two facility location problems under geographical distributed demand where customer' choices do not necessarily follow RUMs. Each problem is characterized by a distinct source of uncertainty. We also show that our framework is able to handle the case where the reward is stochastic.

\item Unified approximate solution framework. Leveraging an easy-to-implement scenario sampling scheme, we derive a deterministic MILP, which is able to tightly approximate the original stochastic model. However, the MILP can be computationally demanding owing to the growth of problem size as the result of  the scenario sampling. To tackle this challenge, we develop a branch-and-cut Benders decomposition algorithm, which decomposes the MILP into two components: a master problem for the planner and a choice subproblem for the customers. We then derive strong and stable Benders cuts in analytical forms for both integer and fractional solutions. Additionally, we introduce a heuristic cut separation designed to transform fractional solutions into feasible integer solutions while simultaneously generating appropriate cuts. Since our solution approach works on the MILP, it serves as an approximate (heuristic) solution to the original stochastic problem. To validate the solution quality, we provide a methodology to estimate the solution gap through statistical inference.

\item Comprehensive experiments and computational insights:  We conduct extensive computational experiments to test the solution quality and efficiency of our framework on featured applications. Our framework consistently delivers high-quality solutions with reasonable sample sizes. Specifically, in applying our framework to the CAOP under the exponomial choice model, we observe significant increases in the expected revenue; moreover, when applied to the CAOP under the mixed multinomial logit model, our framework achieves a good balance between computational time and solution quality. For other applications, our framework demonstrates superior scalability and is able to solve large-scale instances more effectively than current state-of-the-art approaches. Finally, we conduct experiments to shed light on the impact of utility-reward relation on computational difficulty. These results reveal that a stochastic choice-based DPP instance will be more tractable if the correlation between rewards and utilities becomes more positive.
\end{itemize}

The remainder of this paper is organized as follows. Section \ref{sec:prob} gives the problem description and three application problems. Section \ref{sec:SBBD} proposes a sampling-based approximate formulation for the problem, which is then solved by a decomposition approach proposed in Section~\ref{sec:BD}. Section \ref{sec:SAA} discusses the statistical inference and gap estimation methods of our framework. Section \ref{sec:efficiency} provides extensive computational experiments to demonstrate the solution quality and efficiency of our framework. Section \ref{sec:insight} provides computational insights on the impact of the reward-utility relation on formulation tightness and computational difficulty.
Finally, we conclude this paper in Section \ref{sec:conclusion}.

\section{Problem description and applications}\label{sec:prob}

This section presents the description and mathematical formulation, along with three specific applications. Throughout the paper, we represent sets by uppercase letters such as $\ml{A}$ and
vectors and matrices using boldface lowercase letters, such as $\bs{x}=(x_1,x_2,...,x_n)^\top$, where the transpose of $\bs{x}$ is denoted by $\bs{x}^{\top}$. By convention, we denote $\mb{R}$ as the real number field and $\mb{B}:=\{0,1\}$ as the binary space.

\subsection{General setting}\label{sec:gp}
We consider a sequential decision-making process of the planner and the customers. Given a universe of options $\mathcal{J} = \{1,2,...,|\mathcal{J}|\}$, the planner must decide a subset of $\mathcal{J}$ to offer to customers to maximize its expected reward. Let $\bs{x}\in\mathbb{B}^{|\mathcal{J}|}$ be the planner's decision variable, where $x_j = 1$ denotes option $j\in\mathcal{J}$ is offered and $x_j = 0$ otherwise. We assume that the planner's decision is restricted by a set of linear constraints, and the overall decision space is defined as
\begin{equation}
\Omega = \left\{\bs{x}\in\mathbb{B}^{|\mathcal{J}|}: \bs{c} \bs{x}\le \bs{d},  \bs{h} \bs{x} = \bs{g}\right\},
\end{equation}
\JJ{where $\bs{c} \in \mb{R}^{m_1\times |\ml{J}|},\bs{d}\in \mb{R}^{m_1},\bs{h} \in \mb{R}^{m_2\times |\ml{J}|},\bs{g}\in \mb{R}^{m_2}$ are matrices and vectors with deterministic entries.}

The above constraints impose that the planner may only offer a limited number of options to customers. However, we do not impose resource capacities to options, i.e., options have sufficient capacities to accommodate demand. We denote a random vector $\bs{\xi}$ as the uncertainties involved in our problem, which follows known probability distribution $\mb{P}$. The utility of option $j\in\ml{J}$ is the random variable, denoted by $U_j(\bs{\xi})$. 

\JJ{Given $\bs{x}$, we define the \textit{choice set} or \textit{offer set} as $\ml{C}(\bs{x}) = \left\{j\in\mathcal{J}:x_j = 1\right\}$. This is the set of options offered by the planner to its customers. We denote $\bs{U}(\bs{\xi})$ as the random utility vector and $\bs{R}(\xi)$ as the random reward vector that depend on $\bs{\xi}$. 
Additionally, we introduce the random indicator function $\mathbb{I}\{j = \argmax_{k \in \ml{C}(\bs{x})}~U_k(\bs{\xi})\}$, whose output is also a random variable that depends on $\bs{x}$ and $\bs{\xi}$. Specifically, given a fixed planning decision $\bs{x} = \bar{\bs{x}}$ and $\bar{\bs{\xi}}$ as a realization of $\bs{\xi}$, we have $\mathbb{I}\{j = \argmax_{k \in \ml{C}(\bar{\bs{x}})}~U_k(\bar{\bs{\xi}})\} = 1$ if  option $j\in\mathcal{J}$ is the option with the highest utility $U_j(\bar{\bs{\xi}})$ among those offered, and 0 otherwise. If customer selects option $j$, then it will bring a reward of $R_j(\bar{\bs{\xi}})$ to the planner. }

We denote  
\begin{align}
Q(\bs{x},\bs{\xi}) = \sum_{j\in\mathcal{J}}R_j(\bs{\xi}) \cdot \mb{I}\{j = \argmax_{k \in \ml{C}(\bs{x})}~U_k(\bs{\xi})\}.
\end{align} as the random revenue obtained by the planner. 

Therefore, the planner's expected reward over $\bs{\xi}\sim \mb{P}$ under arbitrary decision $\bs{x}$ is
\begin{align}
\E_{\bs{\xi}\sim \mb{P}}[Q(\bs{x},\bs{\xi})] =  \int_{\bs{\xi}} \sum_{j\in\mathcal{J}}R_j(\bs{\xi}) \cdot \mb{I}\{j = \argmax_{k \in \ml{C}(\bs{x})}~U_k(\bs{\xi})\}\cdot p(\bs{\xi})d\bs{\xi},
\end{align}
where $p(\bs{\xi})$ is the probability density function of $\bs{\xi}$, associated with $\mb{P}$.

Altogether, we define the general stochastic choice-based DPP as
\begin{equation}
\textbf{[DPP]}\quad v^* = \max_{\bs{x}\in\Omega}\quad \int_{\bs{\xi}} \sum_{j\in\mathcal{J}}R_j(\bs{\xi}) \cdot \mb{I}\{j = \argmax_{k \in \ml{C}(\bs{x})}~U_k(\bs{\xi})\}\cdot p(\bs{\xi})d\bs{\xi}
\end{equation}

In this paper, we focus on problems with the following characteristics: (i) The planner's decision space is nontrivial, i.e., $|\Omega| \ge 2$; (ii) $\bs{\xi}$ follows continuous distribution; and (iii) $\E_{\bs{\xi}\sim\mb{P}}(\cdot)$ is lower semicontinuous and the reward $R_j(\bs{\xi})$ is bounded, $\forall j\in \ml{J}$ and $\bs{\xi}\sim\mb{P}$. That is, $\underline{R}\le R_j(\bs{\xi})\le \bar{R}$ for real-valued $\underline{R}$ and $\bar{R}$. Given this, we have $\E_{\bs{\xi}\sim\mb{P}}[|Q(\bs{x},\bs{\xi})|] < \infty$ because $\underline{R}|\ml{J}|\le Q(\bs{x},\bs{\xi})\le \bar{R}|\ml{J}|$.

Below, we showcase three existing applications that can be effectively represented by this modeling framework. More specially, the first application is associated with idiosyncratic (internal) uncertainties, and the remaining two are associated with intrinsic (external) uncertainties. However, we emphasize that the model's applicability extends beyond these specific examples.

\subsection{Application: Constrained assortment optimization under random utilities}\label{sec:app}

In our first demonstration, we focus on the constrained assortment optimization problem (CAOP) under random utilities. Consider a retail firm that has to plan the assortment of products to offer to its customers. The retailer has access to a universe of products denoted by the set $\ml{A}$, and can offer at most $\tau$ products. We use set $\ml{S}$ with $|\ml{S}| \leq \tau$ to denote the selected assortment. Each product $a \in \ml{A}$ has a unit profit margin (or selling price\footnote{Without loss of generality, we can assume that the products have a zero cost and therefore $R_a$ can represent the price offered to customers. We refer readers to \citet{leitner2024exact} for CAOP with product cost.}) $R_a$. Naturally, a customer may decide to buy none of the products offered -- this is captured by the outside (no-purchase) option which is denoted as $\{0\}$. Total number of options is $|\ml{A}|$ + 1, and the no-purchase option is always available.

Without loss of generality, we can assume that the utility of alternative $a\in  \ml{A}\cup \{0\}$ is defined as $U_a(\bs{\xi}) = V_a(\bs{\xi_1}) + \bs{\xi_2}$, where $\bs{\xi} = (\bs{\xi}_1, \bs{\xi}_2)$ consists of two random components. Here, $V_a(\bs{\xi_1})$ is referred to as the ideal utility which depends on $\bs{\xi}_1$, and $\boldsymbol{\xi}_2$ reflects customer' idiosyncratic preferences. For any given realization of $\bs{\xi}$, the customer chooses the option that has the highest utility, i.e., the function $\mb{I}\{a = \argmax_{k \in  \ml{S} \cup \{0\}}~U_a(\bs{\xi})\}$ indicates whether option $a$ is chosen. As a result, the generic model for the CAOP can be formulated as
\begin{align}\label{eqt:CAOP}
\textbf{[CAOP]} \quad \max_{\ml{S}\subseteq  \ml{A}, | \ml{S}| \leq \tau}\int_{\bs{\xi}}\sum_{a\in  \ml{S}}R_a \cdot \mb{I}\{a = \argmax_{k \in  \ml{S} \cup \{0\}}~U_a(\bs{\xi})\} \cdot p(\bs{\xi})d\bs{\xi}
\end{align}
where the objective is to maximize the expected revenue through selecting $\ml{S}$ from $\ml{A}$.

To adapt our framework to CAOP, we first redefine set $\ml{J} = \ml{A}\cup \{0\}$. Now, let $x_{j}$ be the binary variable that is equal to 1 if option $j$ exists in the option set, and 0 otherwise. Then, it is easy to see that with the decision space being $\bs{x} \in \Omega:=\left\{\bs{x}\in \B^{|\ml{J}|}:    \sum_{j \in \ml{J}} x_j \le \tau + 1, x_0 = 1\right\}$ and $\ml{C}(\bs{x}) = \{j\in \ml{J}: x_j = 1\}$, the CAOP automatically satisfies our framework in Section \ref{sec:gp}.

Problem~(\ref{eqt:CAOP}) provides the general framework for modeling a broad spectrum of cardinality-constrained assortment optimization problems under RUMs. In principle, the random component $\bs{\xi}$, which causes the randomness in consumers' utilities, can follow any probability distributions. By restricting $\bs{\xi}$ to certain utility distribution families, the resultant CAOP becomes computationally more tractable. For example, with the utility defined as $U_a(\bs{\xi}) = V_a + \bs{\xi_2}$, where $V_a$ is deterministic and $\bs{\xi_2}$ follows i.i.d. Gumbel distribution, we have the CAOP under the multinomial logit model, which is solvable in polynomial time \citep{rusmevichientong2010dynamic}. However, when $\bs{\xi}$ conforms to other distributions, leading to alternative choice models like the exponomial choice model~\citep{alptekinouglu2016exponomial}, the CAOP problem becomes NP-hard, rendering the task of deriving high-quality approximate solutions a significant challenge~\citep{aouad2023exponomial}.
 

\subsection{Application: Discrete location-pricing with geographically distributed demand}\label{sec:app2}

The second application is a facility location and discrete pricing problem (FLoP), which is a stochastic variant and extension of existing works on facility mill pricing~\citep{kochetov2015comparison,lin2023facility,panin2014complexity}. The problem is characterized by geographically distributed demand, where customers are spread across a geographic area. Each customer is associated with a specific location, represented by latitude and longitude coordinates, denoted as $\bs{\xi} = (\text{lon}, \text{lat})$. We model the distribution of customers across the geographic area as a probability distribution, denoted by $\bs{\xi} \sim \mathbb{P}$. Here, $\mathbb{P}$ signifies the probabilistic framework governing the spatial distribution of customers.

Candidate facilities are at discrete sites, denoted by set $\mathcal{A}$. Each facility, once open, imposes a service charge to customers who patronize it. The service charge at each facility is to be selected from the set of discrete (bounded) pricing levels $\ml{L}$,  with $s_l$ being the charge at level $l$, $\forall l \in \ml{L}$.

For a customer located at $\bs{\xi}$, the total cost to seek service from a specific facility $a$ at pricing level $l$  is given by $G_{al}(\bs{\xi}) = \alpha d_a(\bs{\xi}) + s_{l}$, where $d_a(\bs{\xi})$ is the shortest distance to facility $a$, and  $\alpha > 0$ is a parameter that transfers the distance into cost. Each customer has a budget or reservation price denoted by $B$. This reservation price is the maximum amount a customer is willing to spend to obtain the service. In other words, it can also be seen as the cost of the outside option provided by a competing company. The customers' goal is to minimize their costs, either by selecting the open facility with the lowest cost or by opting for the outside option (alternative service) if it offers a lower cost compared with any of the open facility options.

Now, the company needs to open $\tau$ facilities and determine one pricing level for each open facility to maximize its estimated revenue. To this end, we use binary variable $X_{al}$,  $\forall (a,l)\in \mathcal{A}  \times \mathcal{L}$, to denote the location-pricing decision, which is 1, if facility $a$ is open at pricing level $l$ with $s_{l}$ being the corresponding service charge. Then the decision space is $
\left\{\bs{X}\in \B^{|\ml{A}|\times |\ml{L}|}:  \sum_{l \in \ml{L}}X_{al} \leq 1, \forall a \in \ml{A}, \sum_{a \in \ml{A}}\sum_{l \in \ml{L}}X_{al} = \tau   \right\}$,
which imposes that  $\tau$ facilities are open, and each open facility must set a pricing level.

For each location-pricing pair $(a,l) \in \ml{A}  \times \ml{L}$, customers perceive a cost of $ G_{al}(\bs{\xi})$, and the company observes a revenue of $r_{al} = s_{l}$ if the facility is patronized by customers. We define auxiliary parameters $G_{00} = B$ and $r_{00} = 0$ to indicate the cost and reward of the outside option. With these notations, we have that $\ml{C}(\bs{X}) = \left\{(a,l) \in \ml{A}  \times \ml{L}: X_{al} = 1\right\} \cup \{(0,0)\}$. Customers' choices can then be modeled as  $\mb{I}\{(a,l) = \argmax_{(a^\prime,l^\prime) \in\ml{C}(\bs{X}) }-G_{a^\prime l^\prime}(\bs{\xi})\}$. We can  then represent the revenue-maximizing decision problem by
\begin{align}\label{prob:FLoP}
\textbf{[FLoP]} \quad \max_{\bs{X}\in\Omega}\int_{\bs{\xi}} \left(\sum_{a\in \ml{A}}\sum_{l \in \ml{L}}r_{al} \cdot \mathbb{I}\{(a,l) = \argmax_{(a^\prime,l^\prime) \in \ml{C}(\bs{X})}-G_{a^\prime l^\prime}(\bs{\xi})\} \right) p(\bs{\xi})d\bs{\xi}.
\end{align}

Finally, we redefine the parameters, variables, decision sets and decision space as
\begin{align*}
& \bs{U} = \{-G_{00}, -G_{11}, -G_{12},...,-G_{1|\ml{L}|}, -G_{21}, -G_{22},...,-G_{2|\ml{L}|},...,-G_{|\ml{A}||\ml{L}|}\} \\
& \bs{R} = \{r_{00}, r_{11}, r_{12},...,r_{1|\ml{L}|}, r_{21}, r_{22},...,r_{2|\ml{L}|},...,r_{|\ml{A}||\ml{L}|}\} \\
& \bs{x} = \{X_{00}, X_{11}, X_{12},...,X_{1|\ml{L}|}, X_{21}, X_{22},...,X_{2|\ml{L}|},...,X_{|\ml{A}||\ml{L}|} \}\\
& \ml{J}^a = \{a|\ml{L}| + 1, \dots, a|\ml{L}|+|\ml{L}|\}, \forall a\in \tilde{\ml{A}} = \{0,\dots, |\ml{A}| - 1\}\\
& \ml{J} = \cup_{a\in\tilde{\mathcal{A}}}\mathcal{J}^a \cup \{0\}\\
&\ml{C}(\bs{x}) = \{j\in\ml{J}:x_j = 1\}\\
&\Omega =  \left\{\bs{x}\in \B^{|\ml{J}|}: x_{0} = 1, \sum_{j \in\ml{J}} x_{j} = \tau + 1,  \sum_{j\in\ml{J}^a}x_{j} \leq 1, \forall a \in \tilde{\ml{A}}  \right\}
\end{align*}
With them, Problem (\ref{prob:FLoP}) becomes the unified form
\begin{align}
\max_{\bs{x}\in\Omega}\int_{\bs{\xi}}\sum_{j\in\ml{J}}R_j \cdot \mathbb{I}\{j = \argmax_{k \in \ml{C}(\bs{x})}~U_{k}(\bs{\xi})\} \cdot p(\bs{\xi})d\bs{\xi}.
\end{align}

\subsection{Application: Market share maximizing facility location problem}\label{sec:app3}

In the previous two applications, the reward $R_j$ for option $j \in \ml{J}$ is constant and does not depend on any uncertainty parameter. However, in many real-world problems, $R_j$ may exhibit significant dependence on the uncertainty realization. Our third application example addresses one of such problems, known as the Market Share Maximizing Facility Location Problem (MSMFLP) under the partially binary choice rule.

Similar to FLoP, customers in MSMFLP are distributed across a geographical area, and their locations are represented as $\bs{\xi} = (\text{lon}, \text{lat}) \sim \mathbb{P}$. Facilities are located at discrete sites, denoted by set $\mathcal{J}$. A company is to select $\tau$ facilities to open so as to attract as many customers as possible. We use set $\mathcal{S}$ to denote the set of open facilities. Customers' patronizing behaviors to open facilities are interpreted probabilistically and forecasted using the partially binary (logit) rule~\citep{fernandez2017new,mendez2023store}. Specifically,  we use the well-known Huff-like (gravity-based) utility and compute the perceived utility of a customer at $\bs{\xi}$ seeking service from facility $j \in \mathcal{J}$ as $U_j(\bs{\xi}) = q_j d^{-2}_j(\bs{\xi})$, where $q_j>0$ refers to the attractiveness/service quality of facility $j$, and $d_j(\bs{\xi})$ is the distance between the customer's location and facility $j$. Additionally, there is a competitor in the market that offers an equivalent service (i.e., the outside option), whose utility to customers is represented as $O\ge 0$.

Under the partially binary rule, customers will consider the most attractive facility operated by the company and the competing option. The probability of customers using the service from the company (i.e., the market share)  is  $\frac{\max_{j \in \mathcal{S}} U_{j}(\bs{\xi})}{\max_{j \in \mathcal{S}} U_{j}(\bs{\xi}) + O} $, which leads to the decision problem as follow:
\begin{align}
\textbf{[MSMFLP]} \quad \max_{\mathcal{S}\subseteq \ml{J}, |\mathcal{S}|=\tau}\int_{\bs{\xi}} \frac{\max_{j \in \mathcal{S}} U_{j}(\bs{\xi})}{\max_{j \in \mathcal{S}} U_{j}(\bs{\xi}) + O} \cdot p(\bs{\xi})d\bs{\xi}
\end{align}
To represent the problem using our framework, we define the reward as
\begin{align} \label{R_on_xi}
R_{j}(\bs{\xi}) = \frac{U_{j}(\bs{\xi})}{U_{j}(\bs{\xi}) + O},
\end{align}
which integrates the effect of the outside option. Note that since $U_j(\bs{\xi}) > 0$ and $O\ge 0$, the reward is always bounded with $0 < R_j(\bs{\xi}) \le 1$.
With the decision space as $\bs{x} \in \Omega:=\left\{\bs{x}\in \B^{|\mathcal{J}|}: \sum_{j \in \mathcal{J}} x_j = \tau\right\}$, and the set of available options $\ml{C}(\bs{x}) = \{j\in \mathcal{J}: x_j = 1\}$, we can thus restate the problem as
\begin{align}
\max_{\bs{x} \in \Omega} \quad & \int_{\bs{\xi}} \sum_{j \in \mathcal{J}} R_j(\bs{\xi}) \cdot \mb{I}\{j = \arg\max_{k \in \ml{C}(\bs{x})} U_k(\bs{\xi})\} \cdot p(\bs{\xi}) \, d\bs{\xi},
\end{align}
where $R_j(\bs{\xi})$ depends on the uncertainty parameter $\bs{\xi}$ as showed in (\ref{R_on_xi}).

\subsection{Summary of distinct features of three applications}
Apart from sources of uncertainties, the above three applications differ in the relation between reward and utility. In particular, in CAOP, there is no obvious correlation between the reward and the utility as the reward is exogenous; in FLoP, the reward and the utility are negatively related, i.e., increasing the pricing level of a facility reduces its utility; in MSMFLP, there is a strong positive relation between the reward and the utility since a facility's reward is directly proportional to its utility as indicated in (\ref{R_on_xi}). This feature significantly affects the computational difficulty. We will discuss it in details in later sections.

To avoid confusion with different applications, we summarize them in the following table.
\begin{table}[htb]
\centering
\caption{Distinct features of three applications after being represented by our framework}
\label{tab:features}
\resizebox{\textwidth}{!}{%
\begin{tabular}{ccccccc}
\toprule[1.5pt]
Application            &  & CAOP                       &  & FLoP                                                                &  & MSMFLP             \\ \cline{1-1} \cline{3-3} \cline{5-5} \cline{7-7} 
Planner's decision set &  & $\ml{A}$                   &  & $\ml{A}\times\ml{L}$                                                &  & $\ml{J}$           \\
Effective outside option         &  & $\{0\}$                    &  & $\{(0,0)\}$                                                         &  & $\emptyset$        \\
Customer's choice set  &  & $\ml{J} = \ml{A}\cup\{0\}$ &  & $\ml{J} = \cup_{a\in\tilde{\mathcal{A}}}\mathcal{J}^a \cup   \{0\}$ &  & $\ml{J}$           \\
Decision space &
   &
  $\left\{\bs{x}\in\mb{B}^{|\ml{J}|}:x_0 = 1, \sum\limits_{j \in   \ml{J}} x_j \le \tau + 1\right\}$ &
   &
  $\left\{\bs{x}\in\mb{B}^{|\ml{J}|}:x_{0} = 1, \sum\limits_{j   \in\ml{J}} x_{j} = \tau + 1,    \sum\limits_{j\in\ml{J}^a}x_{j} \leq 1, \forall a \in   \tilde{\ml{A}}  \right\}$ &
   &
  $\left\{\bs{x}\in\mb{B}^{|\ml{J}|}:\sum\limits_{j \in   \mathcal{J}} x_j = \tau\right\}$ \\
Uncertain reward       &  & No                         &  & No                                                                  &  & Yes                \\
Reward-utlity relation &  & Not related                &  & Negatively related                                                  &  & Strongly positively related \\ \bottomrule[1.5pt]
\end{tabular}%
}
\end{table}

\section{Unified modeling and solution framework}\label{sec:SBBD}

In \textbf{[DPP]}, deriving an explicit formulation for the expectation and integral is challenging due to the specific probability distributions governing $\bs{\xi}$, rendering general analysis intractable. In this section, we leverage sample average approximation (SAA) to approximate \textbf{[DPP]} with a finite number of scenarios of $\bs{\xi}$. Under SAA, we approximate the uncertain customers' choice behaviors by a deterministic rank list model and formulate an explicit MILP that allows for generating near-optimal solutions to \textbf{[DPP]}. On this basis, we propose a two-stage Benders decomposition approach with customized Benders cuts to efficiently solve the resulting MILP. We name the complete approximation modeling and solution framework as \textit{sampling-based Benders decomposition} (SBBD), details of which will be discussed below. Note that the mathematical proofs for the propositions mentioned in this section are available in Online Appendix A.1.

\subsection{MILP formulation based on SAA}

We sample a finite set of scenarios $\mathcal{I} = \left\{1,2,\dots, N \right\}$ from the distribution of $\bs{\xi}$, with $r_{ij} \in \R^{N\times|\mathcal{J}|}$ and $u_{ij} \in \R^{N\times|\mathcal{J}|}$ be the realization of $R_j(\bs{\xi})$ and $U_j(\bs{\xi})$ under the $i$-th scenario. Hereafter, we interchangeably refer to $N$  as ``number of scenarios'' and ``sample size''.
Under each scenario, each customer has a fixed preference ranking to options according to the realized utilities $\bs{u}$ and chooses the option with the highest ranking (realized utility), leading to a rank list model.

 We denote $\bs{y}$ as the choice variable for the customers, where $\bs{y}$ is related to $\bs{x}$ as
 \begin{equation}
y_{ij} = \begin{cases}
\mathbb{I}\{j = \arg\max_{j^\prime \in \ml{C}(\bs{x})}~u_{ij^\prime}\} &, \text{ if } j\in \ml{C}(\bs{x}) \\
0 &, \text{ otherwise }
\end{cases},\forall i\in \mathcal{I}. \label{set:choice}
\end{equation}
Next, our task is to explicitly formulate the relation between $\bs{x}$ and $\bs{y}$.  If an option is not offered by the planner, then the customers cannot choose it, i.e.,
\begin{equation}
y_{ij}\le x_j,\forall i\in \mathcal{I},\forall j\in \mathcal{J}.\label{constr:yx}
\end{equation}
To ensure that one option will be chosen under each scenario, we impose
\begin{equation}
\sum_{j\in \ml{J}}y_{ij} = 1,\forall i\in \mathcal{I}.
\end{equation}
Then the highest-ranking logic can be modeled as
\begin{equation}\label{constr3}
\sum_{j\in \ml{J}:u_{ij} < u_{ik}}y_{ij} + x_k \le 1, \forall i \in \mathcal{I}, \forall k \in \mathcal{J}.
\end{equation}
It is easy to see that, if an option $k\in \mathcal{J}$ if offered ($x_k = 1$), then $y_{ij} = 0$ for all $j\in\mathcal{J}$ such that $u_{ij} < u_{ik}$; otherwise, if $x_k = 0$, (\ref{constr3}) is inactive.

The choice variable $\bs{y}$ should be binary by definition. However,  since we assume that $\bs{\xi}$ follows continuous distributions, $\bs{y}$ can be relaxed to $ \bs{y} \ge 0$. This is because the realized utilities, after sampling $\bs{\xi}$ from a continuous distribution, are unique almost surely. The fractional solution of $\bs{y}$ can be naturally excluded because the maximum utility option will be unique almost surely.

  Altogether, we can approximate \textbf{[DPP]} by the following MILP formulation
\begin{subequations}\label{model:ADPP}
\begin{alignat}{10}
\label{SABF:obj} v_N =  \max_{\bs{x},\bs{y}}~~ & \frac{1}{N}\sum_{i\in \mathcal{I}}\sum_{j\in \mathcal{J}}r_{ij}y_{ij}\\
\text{s.t.}~~& \bs{x} \in \Omega, \label{SHP:constr0}\\
\textbf{[ADPP]} \qquad \qquad \qquad &\sum_{j\in \ml{J}}y_{ij} = 1,\forall i\in \mathcal{I}, \label{SHP:constr1}\\
&y_{ij}\le x_j,\forall i\in \mathcal{I},\forall j\in \mathcal{J}, \label{SHP:constr2} \\
& \sum_{j\in \ml{J}:u_{ij} < u_{ik}}y_{ij} + x_k \le 1, \forall i \in \mathcal{I},\forall k \in \mathcal{J}, \label{constr:SHP}\\
& y_{ij} \ge 0,\forall i\in \mathcal{I},\forall j\in \mathcal{J}.\label{SHP:constr4}
\end{alignat}
\end{subequations}
Here, \textit{$v_N$ denotes the approximate optimal value of \textbf{[DPP]} under sample size $N$}. 

\subsection{Benders decomposition}
\label{sec:BD}

The linear structure of \textbf{[ADPP]} comes at the cost of increased problem size, particularly when $N$ is large to achieve a high level of approximation quality. As a MILP, \textbf{[ADPP]} can be directly solved by off-the-shelf solvers such as CPLEX and Gurobi. Nevertheless, dealing with large-scale instances may still be challenging.  Benders decomposition emerges as a compelling algorithm in this context, because it allows for projecting out the choice variable $\bs{y}$ and conducting the main solution process only with the planner's variable $\bs{x}$. Specifically, the Benders master problem on the $x$-space  is defined as
\begin{align}
\max \quad &\frac{1}{N} \sum_{i\in \mathcal{I}}\theta_{i}\\
\textbf{[MP]}\qquad\text{s.t.}\quad & \bs{x}\in\Omega,\\
\label{eqt:benders_convex}& \theta_{i} \le \phi_{i}(\bs{x}),\forall i\in \mathcal{I},
\end{align}
where $\phi_i(\bs{x})$ is the expected reward obtained from the choice subproblem for the customers, under scenario $i\in\mathcal{I}$, i.e.,
\begin{subequations}
\begin{alignat}{10}
\phi_{i}(\bs{x}) = \max \quad&\sum_{j\in \mathcal{J}}r_{ij}y_{ij},\\
\text{s.t.}\quad &\sum_{j\in \mathcal{J}}y_{ij} = 1, \label{sp:constr1} \\
\textbf{[SP$_i$]}\qquad \qquad \qquad &y_{ij}\le x_j,\forall j\in \mathcal{J}, \label{sp:constr2}\\
&\sum_{j\in \mathcal{J}:u_{ij} < u_{ik}}y_{ij} + x_k \le 1,\forall k \in \mathcal{J},\label{sp:constr3}\\
&y_{ij}\ge 0,\forall j\in \mathcal{J}.
\end{alignat}
\end{subequations}

The core of the Benders decomposition centers on the generation of linear Benders cut to approximate (\ref{eqt:benders_convex})  and to cut off non-optimal solutions when a planning decision $x$ is given, a process widely known as Benders separation. In the remainder of this section, we describe how we analytically conduct Benders separations at both integer and fractional solutions of $\bs{x}$. Besides, we also introduce a heuristic Benders separation that can further expedite the convergence of Benders decomposition.

\subsubsection{Benders separation with integer solutions}
 Let $\bs{\lambda}\in\mathbb{R}^{N}$, $\bs{\nu}\in\mathbb{R}_{+}^{N\times|\mathcal{J}|}$ and $\bs{\mu}\in\mathbb{R}_{+}^{N\times|\mathcal{J}|}$ be the corresponding dual variables associated with Constraints (\ref{sp:constr1}) - (\ref{sp:constr3}).  The Benders dual subproblem for scenario $i\in\mathcal{I}$ reads
\begin{subequations}
\begin{alignat}{10}
\phi_{i}(\bs{x}) = \min \quad&\lambda_i + \sum_{j\in \mathcal{J}}\mu_{ij} + \sum_{j \in \mathcal{J}} (\nu_{ij}-\mu_{ij}) x_{j},\\
\textbf{[DSP$_i$]}\qquad \qquad \text{s.t.}\quad &r_{ij} - \lambda_i - \nu_{ij} - \sum_{k\in \mathcal{J}:u_{ik} > u_{ij}}\mu_{ik} \le 0,\forall j\in \mathcal{J},\\
& \bs{\nu}_i, \bs{\mu}_i \ge 0
\end{alignat}
\end{subequations}
For a feasible $\bs{\bar{x}}\in\Omega$ , we solve \textbf{[DSP$_i$]} to obtain its optimal solution $(\bar{\lambda}_i, \bar{\bs{\nu}}_i, \bar{\bs{\mu}}_i)$ and generate integer Benders cuts in the form of
\begin{equation}\label{cut:int}
\theta_i \le \bar{\lambda}_i + \sum_{j\in \ml{J}}\bar{\mu}_{ij} + \sum_{j\in\ml{J}}(\bar{\nu}_{ij} - \bar{\mu}_{ij})x_j,\forall i\in \mathcal{I}.
\end{equation}
The following proposition states that  \textbf{[DSP$_i$]} can be solved analytically, requiring only simple sorting and matrix operations; therefore, integer Benders cut (\ref{cut:int}) can be generated efficiently.
\begin{prop} \label{prop:dsp}
For $\bs{\bar{x}}\in\Omega$, define $\mathcal{J}^o = \left\{j\in\mathcal{J}: \bar{x}_j = 1\right\}$ and $j_i^* = \arg\max_{j'\in \mathcal{J}^o}u_{ij'}$. An optimal solution to \textbf{[DSP$_i$]} is given by
\begin{align}
&\bar{\lambda}_{i} = r_{ij_i^*},\\
&\bar{\mu}_{ij} = \begin{cases}
[\max_{j\in \mathcal{J}^o\setminus \left\{j_i^*\right\}   }r_{ij}-\bar{\lambda}_{i}]_{+}  & ,j = j_i^*, \\
  0& ,\forall j\in \mathcal{J}\setminus \left\{j_i^*\right\},
  \end{cases}\\
&\bar{\nu}_{ij} = \begin{cases}
0  & ,\forall j\in \mathcal{J}^o, \\
[r_{ij}-\bar{\lambda}_{i} - \sum_{k\in \mathcal{J}:u_{ik} > u_{ij}}\bar{\mu}_{ik}]_{+}  & , \forall j\in \mathcal{J}\setminus\mathcal{J}^o.
\end{cases}
\end{align}
where $[z]_+ = \max\left\{z,0\right\}$.
\end{prop}

\subsubsection{Benders separation at fractional solutions}

While Benders separation at integer solutions guarantees optimal convergence, introducing Benders cuts at fractional solutions can significantly improve the algorithm's performance. Similar to the previous section, our goal is to derive cuts in analytical forms. However, when $\bs{x}$ is fractional, determining the optimal dual solution via the original formulation \textbf{[DSP$_i$]} proves to be complex and inefficient. To circumvent this challenge, we introduce an alternative, computationally more efficient formulation of the Benders subproblem, i.e.,
\begin{subequations}
\begin{alignat}{10}
\phi_{i}(\bs{x}) = \max \quad&\sum_{j\in \mathcal{J}}r_{ij}y_{ij},\\
\textbf{[nSP$_i$]}\qquad \qquad \quad\text{s.t.}\quad &\sum_{j\in \mathcal{J}}y_{ij} = 1, \label{nsp:constr1} \\
&y_{ij}\le \beta_{ij},\forall j\in \mathcal{J}, \label{nsp:constr2}\\
&y_{ij}\ge 0,\forall j\in \mathcal{J},
\end{alignat}
\end{subequations}
where $\beta_{ij} = \min\left\{x_j, \min_{k\in \mathcal{J}}\left\{1 - \delta_{ijk}x_k\right\} \right\},\forall i\in \mathcal{I},\forall j\in \mathcal{J}$ and \begin{equation}
\delta_{ijk} =
 \begin{cases}
1& \text{if}~ u_{ik} >  u_{ij} \\
0 & \text{otherwise}
\end{cases}, \forall i \in \mathcal{I}, j \in \mathcal{J}, k \in \mathcal{J}
\end{equation} is the ranking parameter.

\vspace{2mm}
\begin{prop}\label{prop:newCAC}
\textbf{[nSP$_i$]} is equivalent to \textbf{[SP$_i$]} for all $\bs{x}\in\Omega$.
\end{prop}
\vspace{2mm}

 The equivalence of \textbf{[nSP$_i$]}  and \textbf{[SP$_i$]} holds when $\bs{x}$ is integer and feasible. Obviously, \textbf{[nSP$_i$]} is a \textit{knapsack} problem where $r_{ij}$ is the value of each item and $\beta_{ij}$ is the weight, subject to total capacity of 1. We can leverage this specific structure of the new subproblem to quickly generate Benders cuts. Let $\bs{\lambda'} \in \mathbb{R}^{N}$ and $\bs{\eta} \in \mathbb{R}_{+}^{N\times |\mathcal{J}|}$ be the dual variables associated with Constraints (\ref{nsp:constr1}) and (\ref{nsp:constr2}). The new dual subproblem is
\begin{subequations}
\begin{alignat}{10}
 \min\quad &\lambda'_i + \sum_{j\in \mathcal{J}}\eta_{ij}\beta_{ij}\\
 \textbf{[nDSP$_i$]}\qquad\text{s.t.}\quad & r_{ij} - \lambda'_i-\eta_{ij} \le 0,\forall j\in \mathcal{J},\\
& \eta_{ij} \ge 0,\forall j\in\mathcal{J}.
\end{alignat}
\end{subequations}
Now we look at fractional $\bs{\bar{x}}$ obtained from the continuous relaxation of \textbf{[MP]}, i.e., $0\le \bar{x}_j\le1,\forall j\in \mathcal{J}$. Let $\bs{\bar{\beta}}$ be the corresponding weight given $\bs{\bar{x}}$. The following result arises.

\vspace{2mm}
\begin{prop}\label{prop:frac1}
Under scenario $i\in I$, sort options in the descending order $\bs{\pi}$ with respect to $\bs{r}$ such that $r_{i\pi_{0}} > r_{i\pi_{1}} > \cdots > r_{i\pi_{|\mathcal{J}| - 1}}$.
Let $j_i^*$ be the critical item such that $\sum_{j' = 1}^{j_i^* - 1}\bar{\beta}_{i\pi_{j'}} < 1 = \sum_{j' = 1}^{j_i^*}\bar{\beta}_{i\pi_{j'}}$.
An optimal dual solution ($\bar{\lambda}'_i,  \bar{\bs{\eta}}_i$) to  \textbf{[nDSP$_i$]} is given by
\begin{align}
&\bar{\lambda}'_i = r_{ij_i^*}\\
&\bar{\eta}_{ij} = [r_{ij} - r_{ij_i^*}]_{+},\forall j\in \mathcal{J}.
\end{align}
\end{prop}
\vspace{2mm}

On this basis, we can rewrite Benders cuts with respect to $\bs{\beta}$ as
\begin{equation}
\theta_i \le \bar{\lambda}'_i + \sum_{j\in \mathcal{J}}\bar{\eta}_{ij}\beta_{ij},\forall i\in \mathcal{I}. \label{BD:beta}
\end{equation}
Note that Cut (\ref{BD:beta}) is implicitly related to $\bs{x}$ by $\bs{\beta}$. According to the definition of $\bs{\beta}$, we can use the \textit{chain rule} to transform Cut (\ref{BD:beta}) into the Bender cut with respect to $\bs{x}$, supported by the following result.

\vspace{2mm}
\begin{prop}\label{prop:frac2}
\JJ{Under each scenario $i\in \mathcal{I}$, given $\bar{\bs{x}}$ and $\bar{\bs{\beta}}$, define
\begin{align}
&\bar{\nu}^\prime_{ij}=\left\{\begin{matrix}
\bar{\eta}_{ij}  ,&\beta_{ij} = \bar{x}_j \\
0  ,& \text{otherwise}
\end{matrix}\right.  \quad \forall i\in \ml{I},\forall j\in \ml{J},\\
&\bar{\mu}^\prime_{ijk} = \left\{\begin{matrix}
 \bar{\eta}_{ij} ,& \bar{\beta}_{ij} = 1-\delta_{ijk}\bar{x}_{k}, k = \arg\min_{k'\in \mathcal{J}} 1-\delta_{ijk'}\bar{x}_{k'}\\
 0, &\text{otherwise}
\end{matrix}\right. \quad \forall i\in \ml{I},\forall j\in \ml{J},\forall k\in \ml{J}.
\end{align}
}
Then Cut (\ref{BD:beta}) is equivalent to
\begin{equation}
\theta_i \le \bar{\lambda}_i' + \sum_{j\in \mathcal{J}}\sum_{k\in\mathcal{J}}\bar{\mu}'_{ijk} + \sum_{j\in \mathcal{J}}(\bar{\nu}'_{ij} - \sum_{k\in J}\delta_{ijk}\bar{\mu}'_{ijk})x_j. \label{cut:frac}
\end{equation}
\end{prop}
\vspace{2mm}

Propositions \ref{prop:frac1} and \ref{prop:frac2} give an easy way to construct valid Benders cuts at fractional solutions based on the new subproblem \textbf{[nSP$_i$]}. Based on Proposition \ref{prop:newCAC}, the derived Benders cuts are also valid at fractional solutions for the original formulation \textbf{[ADPP]}.

\subsubsection{Heuristic separation}
\label{sec:Heuristic_separation}

In addition, we incorporate a heuristic cut separation procedure to expedite the convergence of our Benders approach.  The fundamental idea is to employ a primal heuristic to transform fractional solutions $\bs{\bar{x}}$ into feasible integer solutions $\bs{\tilde{x}}$, following which we activate the integer Benders separation. 
Our primal heuristic is done in a \textit{greedy} manner. Below, we illustrate how to transfer a fractional solution to an integer solution while satisfying constraints.

\begin{itemize}
\item CAOP. We have $x_0 = 1$, $\sum_{j\in \mathcal{J}} x_j  \le \tau + 1$. Given a fractional solution $\bs{\bar{x}}$, we sort $\bs{\bar{x}}$ in descending order to obtain $\bs{\tilde{x}}$. We then proceed to select the largest $\tau + 1$ elements from $\bs{\tilde{x}}$ and set them to 1 (note that the outside option $\{0\}$ is definitely one of them since  $\bar{x}_0 = 1$), while assigning 0 to the remaining elements.
This process enforces a cardinality constraint, where the sum of the selected elements in $\bs{\tilde{x}}$ is equal to the specified value $\tau + 1$ while the outside option $\{0\}$ is always offered.
\eat{
\item Budget constraint $\sum_{j\in \mathcal{J}} c_jx_j  \leq b$. Given a fraction solution $\bar{x}$, we sort it in descending order to obtain $\tilde{x}$. Subsequently, we fill the elements in $\tilde{x}$ with 1 in a greedy manner, starting from the largest values, until the budget is fully utilized. Any remaining elements are set to 0. This approach ensures that the total cost, represented by the weighted sum of elements in $\tilde{x}$, remains within the specified budget limit $b$.
}

\item FLoP. We have $x_{0} = 1, \sum_{j \in\ml{J}} x_{j} = \tau + 1,  \sum_{j\in\ml{J}^a}x_{j} \leq 1, \forall a \in \tilde{\ml{A}}$.
\JJ{Given a fractional solution $\bs{\bar{x}}$, we define two vectors $T, K \in \mb{R}^{|\tilde{\ml{A}}|}$ (both indexed by $a\in \tilde{\ml{A}}$) such that $T_a = \max_{j\in\mathcal{J}^a}\bar{x}_{j}, K_a = \mathop{\arg\max}_{j\in\mathcal{J}^a}\bar{x}_{j}, \forall a\in\tilde{\ml{A}}$.}We then sort the value vector $T$ in descending order to obtain the first $\tau$ elements $t_1,\dots, t_{\tau}$. For each such element, we find the corresponding index $k_1, \dots, k_{\tau}$ from the index vector $K$.
We obtain $\bs{\tilde{x}}$ by setting $\tilde{x}_{k_1} = \dots = \tilde{x}_{k_{\tau}} = 1$ and the remaining elements to 0. This process ensures that each facility $j$ should be associated with at most one price and the overall number of open facilities should be equal to $\tau$.

\item MSMFLP. We have $\sum_{j\in \ml{J}}x_j = \tau$. There is no outside option. Given a fractional solution $\bs{\bar{x}}$, we sort it in descending order to obtain $\bs{\tilde{x}}$. We then proceed to select the largest $\tau$ elements from $\bs{\tilde{x}}$ and set them to 1, while assigning 0 to the remaining elements. This process ensures that the number of open facilities is equal to $\tau$.

\end{itemize}

\subsubsection{Implementation}\label{sec:implement}
Finally, we summarize the implementation process of our two-stage Benders decomposition. We use an iterative process to implement Stage 1 and a branch-and-cut process to implement Stage 2. The heuristic cut separation procedure is only triggered in Stage 2.
\begin{itemize}

\item Stage 1. We first solve the continuous relaxation of the Benders formulation \textbf{[MP]} + \textbf{[nSP$_i$]}, where $\bs{x}\in\Omega':= \left\{\bs{x}\in[0,1]^{|\mathcal{J}|}:\bs{c} \bs{x} \le \bs{d}, \bs{h}\bs{x}=\bs{g} \right\}$ by iterative Benders algorithm and store Benders cuts (\ref{cut:frac}) generated during the process. We adopt an iterative process to implement Stage 1 and stop the iteration process when $\left|UB - LB\right| / LB \le \rho$, where $\rho$ is a small tolerance. We set $\rho  = 10^{-2}$ for CAOP and $\rho  = 10^{-4}$ for FLoP and MSMFLP.

Optionally, we apply a stabilization procedure to mitigate the zig-zag of $x$ solutions to reduce the number of iterations required. Specifically, we introduce a stabilizer $x^s$, which lies in the relative interior space of the convex hull of $\Omega$. At $t$-th iteration, we stabilize the current solution $x^t$ by setting
\begin{align}
&x^t = \varrho^{t} x^t + (1 - \varrho^{t}) x^s,\\
&x^s = (x^s + x^t) / 2,
\end{align}
where $\varrho^{t} \in (0,1]$ is the step size. The use of stabilization generally leads to faster convergence in Stage 1 but does not necessarily improve the solution efficiency for all problem classes. Therefore, unless otherwise specified, the default value of $\varrho^{t} $ is set to 1.

\item Stage 2. We embed Stage 2 in the branch-and-cut process of Gurobi, based on the Benders formulation \textbf{[MP]} + \textbf{[SP$_i$]}. We add the Benders cuts (\ref{cut:frac}) generated in Stage 1 to warm-start the B\&C process. We set \textit{PreSolve} to 0 to prevent Gurobi from removing these Benders cuts. In the branch-and-cut process, we generate Cut (\ref{cut:int}) for all explored integer nodes as lazy cuts within Gurobi's \textit{callback} function. The heuristic separation procedure is carried out for each cut pass at the root node of the search tree and then repeated once every 200 nodes in the subsequent branches. This periodicity helps prevent an excessive number of cuts from being inserted into the problem, maintaining the balance between improving the solution and computational efficiency. Cuts are added only if they violate the current solution $(\bs{\bar{x}},\bs{\bar{\theta}})$ by a minimum relative violation (MRV) which is set as $10^{-5}$ in our implementation. We set \textit{LazyConstraints} to 1 to allow lazy cuts.
\end{itemize}

\section{Estimation of solution gap}\label{sec:SAA}
Our framework utilizes the approximate formulation \textbf{[ADPP]}, which indeed generates a heuristic solution that may be suboptimal. In most cases, the true optimal solution of \textbf{[DPP]} is not known, and it is hard to validate the optimality of a given solution due to the stochastic nature of \textbf{[DPP]}. This section thus discusses estimation methods so that we can gauge the quality of  solutions by our  framework.

Let $\ml{X}^*, \ml{X}_N \subseteq\Omega$ be the set of optimal solutions for \textbf{[DPP]} and \textbf{[ADPP]} with sample size $N$,  respectively.
 According to~\cite{kleywegt2002sample}, the probability of the event $\{\ml{X}_N\subset \ml{X}^*\}$ approaches 1 exponentially fast as $N \rightarrow \infty$ \footnote{To achieve this result, a regularity assumption should be satisfied. To illustrate this, consider a mapping $g:\Omega\setminus\mathcal{X}^*\mapsto \Omega$. We define $H(\bs{x},\bs{\xi}) = Q(g(\bs{x}),\bs{\xi}) - Q(\bs{x},\bs{\xi})$ for all $\bs{x}\in \Omega\setminus\mathcal{X}^*$. Then, $H(\bs{x},\bs{\xi})$ is a bounded random variable, i.e., $(\underline{R} - \bar{R})|\ml{J}|\le H(\bs{x},\bs{\xi}) \le (\bar{R} - \underline{R})|\ml{J}|$, which satisfies \textit{Assumption} (A) in \citet{kleywegt2002sample}}.

However, for large-scale problems, using a large sample size $N$ could result in computationally challenging approximate problems. It is thus necessary to discuss techniques that can estimate the optimality gap of $\bs{x_N}$ for a finite (non-large) $N$ so that we can assess how close $\bs{x}_N$ is to the true optimal solution and validate the effectiveness of the sampling-based framework, especially when dealing with practical problems where using a large $N$ may not be computationally tractable, and it is essential to ensure reliable performance even with a relatively small sample size.

Recall that $v_N$ and $v^*$ denote the optimal values of the problem with sample size $N$ and the true problem.  It can be verified that \footnote{\JJ{For a quick look of $\E[v_N] \geq v^*$, note that for any given $\bs{\bar{x}}$,  it holds that $\E_{\bs{\xi}^i\sim\mb{P}}[Q(\bs{\bar{x}},\bs{\xi}^i)] = \E_{\bs{\xi}\sim\mb{P}} [Q(\bs{\bar{x}},\bs{\xi})]$ since $\bs{\xi}^i$ and $\bs{\xi}$ follow the same distribution. Then, $\max_{\bs{x} \in \Omega} \E \left[\frac{1}{N}\sum_{i=1}^{N}Q(\bs{x},\bs{\xi}^i) \right] = \max_{\bs{x} \in \Omega} \E_{\bs{\xi}\sim\mb{P}}[ Q(\bs{x}, \bs{\xi})] = v^*$. Since $\max_{\bs{x} \in \Omega} \frac{1}{N}\sum_{i=1}^{N}Q(\bs{x},\bs{\xi}^i) \ge \frac{1}{N}\sum_{i=1}^{N}Q(\bs{x},\bs{\xi}^i), \forall \bs{x}\in\Omega$, we have $\E [v_N] = \E \left[\max_{\bs{x}\in\Omega}\frac{1}{N}\sum_{i=1}^{N}Q(\bs{x},\bs{\xi}^i) \right] \ge \max_{\bs{x} \in \Omega} \E \left[\frac{1}{N}\sum_{i=1}^{N}Q(\bs{x},\bs{\xi}^i) \right] = v^*$} }
$\E[v_N] \geq v^*$, leveraging the result in~\cite{shapiro2021lectures}. Therefore, there is a positive expected bias of using the sampling based formulation to approximate the problem, and a statistical upper bound to $v^*$ can thus be estimated by $\E[v_N]$. To estimate $\E[v_N]$, we use a statistical approach: by generating $M$ independent samples of $\bs{\xi}$, each of size $N$, and solving the corresponding [\textbf{ADPP}] for each sample, we can obtain $M$ optimal objective values denoted as $v_N^m$ for $m = 1, 2, ..., M$, along with their corresponding solutions $\bs{x}^m_N$ for $m = 1, 2, ..., M$. Then, the quantity
\begin{align}
\bar{v}^{M}_N = \frac{1}{M} \sum^{M}_{m =1}v_N^{m}
\end{align}
is an unbiased estimator \footnote{\JJ{$\E[\bar{v}^{M}_N] = \frac{1}{M} \sum^{M}_{m =1}\E[v_N^{m}] = \E[v_N]$}} of $\E[v_N]$, which serves as a statistical upper bound to $v^*$. Note that if the samples are i.i.d, then $\E[v_{N+1}] \le \E[v_N]$ \footnote{\JJ{This is because $\E[v_{N+1}] = \E[\max_{\bs{x}\in\Omega}\frac{1}{N+1}\sum_{i=1}^{N + 1}Q(\bs{x}, \bs{\xi}^i)] = \E[\max_{\bs{x}\in\Omega}\frac{1}{N+1}\sum_{i=1}^{N + 1}\frac{1}{N}\sum_{j=1,j\neq i}^{N + 1}Q(\bs{x}, \bs{\xi}^j)] \le \frac{1}{N+1}\sum_{i=1}^{N + 1}\E[\max_{\bs{x}\in\Omega}\frac{1}{N}\sum_{j=1,j\neq i}^{N + 1}Q(\bs{x}, \bs{\xi}^j)] = \E[v_N]$.}}, and thus, the upper bound could be tighter if we increase $N$.

 The variance of the above estimator is
\begin{align}
S_{\bar{v}^{M}_N }^2 = \frac{1}{M(M-1)}\sum^{M}_{m =1}(v_N^{m}  - \bar{v}^{M}_N)^2
\end{align}
which measures the variability of the estimator $\bar{v}^{M}_N$ and can be reduced by increasing $M$.

We can directly compute the lower bound $\hat{v} = \max_{\bs{x}^m_N}\E_{\bs{\xi}\sim\mb{P}}[Q(\bs{x}^m_N, \bs{\xi})]$ if we know the explicit form because each $\bs{x}^m_N$ is a feasible solution to \textbf{[DPP]}. However, in most cases the closed form is not available.
Instead, we use a large sample size $N^\prime$ with a set of scenarios $\{\bs{\xi}^1, \bs{\xi}^2, ..., \bs{\xi}^{N^\prime}\}$ to estimate $\hat{v}$. Then, since any $\bs{x}^m_N$ is a feasible solution, $\frac{1}{N^\prime}\sum_{k =1}^{N^\prime}Q(\bs{x}^m_N,\bs{\xi}^k)$ is naturally a statistical lower bound on $v^*$.  Clearly, under $M$ independent [\textbf{ADPP}], the best estimated solution $\bs{\hat{x}}$ and the best estimated lower bound $\hat{v}$ are
\begin{align}
&\bs{\hat{x}} =  \argmax_{\bs{x}^m_N}~\frac{1}{N^\prime}\sum_{k =1}^{N^\prime}Q(\bs{x}^{m}_N,\bs{\xi}^k)\\
&\hat{v} =\frac{1}{N^\prime}\sum_{k =1}^{N^\prime}Q(\bs{\hat{x}},\bs{\xi}^k)
\end{align}

The estimated variance is
\begin{align}\label{s2v}
&S_{\hat{v}}^2 = \frac{1}{N^\prime(N^\prime-1)} \sum_{k =1}^{N^\prime}(Q(\bs{\hat{x}},\bs{\xi}^k)  -\hat{v})^2
\end{align}
As a result, the total variance of the above estimation process is
\begin{align}
\sigma^2 = S_{\bar{v}^{M}_N }^2 + S_{\hat{v}}^2
\end{align}
which represents the variability of the estimated lower and upper bounds. Lower values of  $\sigma^2$ indicate more consistent and reliable results. For practical usage, $N^\prime$ should be sufficiently large to ensure that the uncertainty is well represented, and the obtained statistical lower bound $\hat{v}$ is reliable (i.e., with small enough $S_{\hat{v}}^2$). Fortunately, it is computationally fast to evaluate $\hat{v}$, making the use of large sample sizes practical. Therefore, we set $N^\prime = 1 \times 10^6$ throughout this paper.  According to (\ref{s2v}), the resulting variance $S_{\hat{v}}^2$ is negligible\footnote{In our experiments, we use $N' = 10^6$. As a result, we observe $S_{\hat{v}}^2 < 10^{-10}$ for all the tested instance.}.
 A a result,  $\sigma^2$ mainly depends on  $S_{\bar{v}^{M}_N }^2$, which could be a non-trivial value since $M$ is small in most of the implementation.

Eventually, we estimate the relative optimality gap as
\begin{align}
\Delta = \frac{\bar{v}^{M}_N - \hat{v}}{\hat{v}} \times 100\%
\end{align}
It is important to note that $\Delta$ is a stochastic term with $\E[\Delta] \geq 0$. In contrast to optimality gaps in deterministic optimization problems, $\Delta$ can possibly take negative values, especially  in situations where a sample size $N$ leads to a near-optimal approximation of the problem, and the number of independent samples $M$ is relatively small (and thus shows high variability). Nevertheless, $\Delta$ still provides a good indication of  the quality of the approximation.  Moreover, since the estimation is characterized by uncertainty and variance, we can construct a $\alpha\%$ unilateral confident upper bound of the relative optimality gap \citep{taherkhani2020benders}, i.e., 
\begin{align}
\Delta_{\alpha} = \frac{\bar{v}^{M}_N - \hat{v} +  z_{1-\alpha} \cdot \sigma}{\hat{v}} \times 100\% = \Delta +  z_{1-\alpha} \cdot \sigma/\hat{v}\times 100\%
\end{align}
where $z_{1-\alpha}$ represents the z-score of the standard normal distribution corresponding to the $(1-\alpha)$th percentile. $\Delta_{\alpha}$ effectively provides a range within which we expect the true optimality gap to lie with a confidence level of $\alpha$.

\remark We use $\Delta$ and $\Delta_{\alpha}$ as metrics to evaluate the solution quality and stability of approximation under $N$, both of which are calculated with $M$ independent replications (sampling), meaning that they are determined after solving $M$ instances of \textbf{[ADPP]}. These two metrics are independent of the true optimal solution which is not known in most cases.
When constraints exist regarding the timing of decisions, we can use one replication, since in practice, conducting one replication with a carefully chosen sample size is sufficient to obtain near-optimal solutions~\citep{kleywegt2002sample}.

\remark  To minimize the variance, one can employ \textit{Latin Hypercube Sampling} for uncertainty sampling, as opposed to \textit{Monte Carlo Sampling}. Our preliminary tests indicate that employing Latin Hypercube Sampling significantly improves approximation accuracy and reduces variance when compared to Monte Carlo Sampling (for details, refer to Online Appendix A.2). Moving forward, the default sampling approach will be Latin Hypercube Sampling.

\section{Computational experiments}\label{sec:efficiency}

In this section, we conduct extensive computational experiments to demonstrate the effectiveness of our proposed SBBD framework in terms of solution quality and computational speed. We divide this section into two main parts: CAOP-related problems and facility location problems. 

Note that the approximate formulation \textbf{[ADPP]} is a challenging optimization problem. To demonstrate the efficiency of our SBBD, we compare the computational speed of our SBBD with other exact solution approaches. Here, we briefly introduce three benchmark solution approaches: \JJ{(i) Directly solving \textbf{[ADPP]} without decomposition (\emph{MILP}).} (ii) Bilevel branch-and-cut approach (\emph{BIBC}). In fact, \textbf{[ADPP]} is a Stackelberg game between the planner and the customers. It can be restated as a bilevel program where the upper-level problem is for the planner to maximize the expected reward, and the lower-level problem is for the customers to maximize their own utilities. This bilevel formulation can be solved by a bilevel branch-and-cut approach~\citep{lin2023facility}. Note that both MILP and BIBC are applicable to all three applications. (iii) Generalized Benders decomposition approach (\emph{GBD}). This approach is designed for MSMFLP, where the problem can be reformulated as a mixed-integer linear fractional program \citep{lin2021branch}.

For performance analysis, we define the following metrics: t[s] refers to the computational time in seconds; \#Node refers to the number of explored nodes in the branch-and-cut tree; \#C refers to the number of cuts added; rgap[\%] refers to the root gap, which reflects the tightness of relaxation bound right after processing the root nodes. Specifically, let $z_{rel}$ be the continuous relaxation bound after processing the root nodes and $z_{opt}$ be the optimal solution, we compute the root gap as $(z_{rel} - z_{opt})/z_{rel} \times 100\%$; ogap[\%] refers to the gap between best upper bound and best lower bound. An instance is said to be solved \textit{optimally} if ogap$<0.01\%$.

All experiments are coded in Python 3.8, with Gurobi 10.0.3 as the solver, and run on a 16 GB Windows PC with an Intel i5-10500 CPU at 3.10 GHz. Throughout this section, the time limit for solving an instance is set as 1 hour.

\subsection{Experiment on CAOP}

For CAOP, we consider two RUMs: \textit{exponomial choice model} and \textit{mixed multinomial logit model}.

\subsubsection{CAOP under exponomial choice model}

 We first apply our SBBD to the CAOP under the exponomial choice model, which has gained increasing attention in recent years~\citep{alptekinouglu2016exponomial, aouad2023exponomial, berbeglia2022comparative}.  In this model, an alternative $j \in \ml{J}$ has a utility $U_j = V_j - \xi_j$, comprising a deterministic component $V_j$ and an an i.i.d. exponential random term with rate $\zeta$, i.e., $\xi_j \sim exp(\zeta)$.
 
 Under the exponomial choice model, the objective function $\E_{\bs{\xi}\sim\mb{P}}[Q(\bs{x},\bs{\xi})]$ of \textbf{[DPP]} has a complex closed-form formulation.  Currently, the state-of-the-art solution approach for the resultant CAOP is the (heuristic) dynamic programming method by \citet{aouad2023exponomial}. Therefore, we compare the performance of our SBBD with  \citet{aouad2023exponomial} algorithm. 

We use the original data generation and algorithm code shared by \citet{aouad2023exponomial}. The tested instances are generated with $|\ml{A}| \in \{9, 19, 29, 39\}$  products plus one outside option. Let $\mathcal{J} = \ml{A}\cup \{0\}$, we have $|\mathcal{J}|\in\{10,20,30, 40\}$ options. The cardinality limit is defined as $\tau = \gamma|\mathcal{J}| - 1$, and the decision space is $\Omega:= \left\{x\in\mb{B}^{|\ml{J}|}: x_0 = 1, \sum_{j\in\mathcal{J}}x_j\le \gamma |\mathcal{J}|\right\}$. For a product $a \in \ml{A}$, the reward $R_a$ is generated from a lognormal distribution with scale $\sigma^R$, and the deterministic utility $V_a$ is generated from a normal distribution with mean 1 and variance $\sigma^U$. For the outside option, we set $R_0$ and $V_0$ to 0. The utilities are then sorted in descending order. The rate parameter $\zeta$ is set as 1. In our experiment, we fix $\gamma = 0.3$ because we find a large $\gamma$ would make the cardinality constraint redundant.

For \citet{aouad2023exponomial} algorithm, we set $\hat{\epsilon} = 0.1$, $\hat{\theta} = 2^{-7}$, which corresponds to the highest precision level tested in \citet{aouad2023exponomial}. Note that the approximation guarantee of  the algorithm may not hold because $\hat{\theta}$ value is not small enough (should be $\frac{\hat{\epsilon}}{2|\mathcal{J}|^2}$ in theory). Setting $\hat{\theta} = 2^{-7}$ is to avoid lengthy computational time. One can check Electronic Companion EC.5 of \citet{aouad2023exponomial} for implementation details and explanation. Here we directly use their original code and follow the same setting.

We use two random seeds in data generation: the instance seed  $S_1$ for generating deterministic utility $\bs{V}$ and reward $\bs{R}$, and the scenario seed $S_2$ for generating the relaxation of the random utility $\bs{\xi}$. For each parameter setting, we change the instance seed $S_1$ to generate 50 different instances for $|\ml{J}| \in \{10,20,30\}$ and 20 different instances for $|\ml{J}| = 40$. For fair comparison, we draw $N = 300$ samples by fixing the scenario seed $S_2 = 1$. This is to ensure that our SBBD yields a unique result for each instance. We then compare the solution of \citet{aouad2023exponomial} algorithm and SBBD with single-run implementation. 

Let $\bs{x^A}$ and $\bs{x^S}$ be the solutions found by \citet{aouad2023exponomial} algorithm and by SBBD, respectively. We then compute the corresponding objective $v^A$ and $v^S$ by the closed-form choice probability of the exponomial choice model. We use $\Delta^S = \frac{v^{ub} - v^S}{v^{ub}}\times 100\%$ and $\Delta^A = \frac{v^{ub} - v^A}{v^{ub}}\times 100\%$ to measure the solution quality for two approaches, where $v^{ub}$ is an upper bound for CAOP. For $|\ml{J}| \in \{10,20\}$, we enumerate all possible solutions to obtain the optimal objective $v^{ub} = v^*$ and set it as the upper bound; whereas for $|\ml{J}| \in  \{30,40\}$, complete enumeration becomes computationally prohibitive. Instead, we estimate the upper bound using the methodology described in Section \ref{sec:SAA}, by varying the scenario seed $S_2$ from 1 to 50, i.e., we have $M = 50$ and set $v^{ub} = \bar{v}_M^N$. Smaller values of  $\Delta^A$ and $\Delta^S$ indicates better solution quality.

\begin{table}[htb]
\centering
\caption{Average $\Delta^A$ and $\Delta^S$ in percentage. `N.A.' indicates that \citet{aouad2023exponomial} algorithm fails to terminate within 1 hour and no solution is returned}
\label{tab:compare-CAOP}
\resizebox{0.6\textwidth}{!}{%
\begin{tabular}{cclcclcclcclcc}
\toprule[1.5pt]
\multicolumn{2}{c}{$|\ml{J}|$} &
   &
  \multicolumn{2}{c}{10} &
   &
  \multicolumn{2}{c}{20} &
   &
  \multicolumn{2}{c}{30} &
   &
  \multicolumn{2}{c}{40} \\ \cline{4-5} \cline{7-8} \cline{10-11} \cline{13-14} 
\multicolumn{2}{c}{$\sigma^R$} &
   &
  0.2 &
  0.5 &
   &
  0.2 &
  0.5 &
   &
  0.2 &
  0.5 &
   &
  0.2 &
  0.5 \\ \cline{1-2} \cline{4-5} \cline{7-8} \cline{10-11} \cline{13-14} 
$\sigma^U=1$ &
  $\Delta^A$ &
   &
  0.00 &
  0.00 &
   &
  0.00 &
  0.01 &
   &
  0.18 &
  0.11 &
   &
  N.A. &
  N.A. \\
 &
  $\Delta^S$ &
   &
  0.06 &
  0.02 &
   &
  0.11 &
  0.06 &
   &
  0.39 &
  0.20 &
   &
  0.37 &
  0.34 \\
\multicolumn{1}{l}{} &
  \multicolumn{1}{l}{} &
   &
  \multicolumn{1}{l}{} &
  \multicolumn{1}{l}{} &
   &
  \multicolumn{1}{l}{} &
  \multicolumn{1}{l}{} &
   &
  \multicolumn{1}{l}{} &
  \multicolumn{1}{l}{} &
   &
  \multicolumn{1}{l}{} &
  \multicolumn{1}{l}{} \\
$\sigma^U=2$ &
  $\Delta^A$ &
   &
  0.00 &
  0.00 &
   &
  0.01 &
  0.00 &
   &
  0.59 &
  0.11 &
   &
  0.46 &
  N.A. \\
 &
  $\Delta^S$ &
   &
  0.05 &
  0.03 &
   &
  0.08 &
  0.10 &
   &
  0.29 &
  0.21 &
   &
  0.34 &
  0.33 \\
\multicolumn{1}{l}{} &
  \multicolumn{1}{l}{} &
   &
  \multicolumn{1}{l}{} &
  \multicolumn{1}{l}{} &
   &
  \multicolumn{1}{l}{} &
  \multicolumn{1}{l}{} &
   &
  \multicolumn{1}{l}{} &
  \multicolumn{1}{l}{} &
   &
  \multicolumn{1}{l}{} &
  \multicolumn{1}{l}{} \\
$\sigma^U=3$ &
  $\Delta^A$ &
   &
  0.00 &
  0.00 &
   &
  0.09 &
  0.06 &
   &
  1.70 &
  2.52 &
   &
  3.53 &
  4.17 \\
 &
  $\Delta^S$ &
   &
  0.04 &
  0.03 &
   &
  0.08 &
  0.08 &
   &
  0.19 &
  0.23 &
   &
  0.29 &
  0.31 \\
\multicolumn{1}{l}{} &
  \multicolumn{1}{l}{} &
   &
  \multicolumn{1}{l}{} &
  \multicolumn{1}{l}{} &
   &
  \multicolumn{1}{l}{} &
  \multicolumn{1}{l}{} &
   &
  \multicolumn{1}{l}{} &
  \multicolumn{1}{l}{} &
   &
  \multicolumn{1}{l}{} &
  \multicolumn{1}{l}{} \\
$\sigma^U=4$ &
  $\Delta^A$ &
   &
  0.71 &
  1.28 &
   &
  1.76 &
  2.84 &
   &
  5.98 &
  6.32 &
   &
  25.43 &
  17.22 \\
 &
  $\Delta^S$ &
   &
  0.04 &
  0.05 &
   &
  0.10 &
  0.06 &
   &
  0.15 &
  0.13 &
   &
  0.27 &
  0.21 \\
\multicolumn{1}{l}{} &
  \multicolumn{1}{l}{} &
   &
  \multicolumn{1}{l}{} &
  \multicolumn{1}{l}{} &
   &
  \multicolumn{1}{l}{} &
  \multicolumn{1}{l}{} &
   &
  \multicolumn{1}{l}{} &
  \multicolumn{1}{l}{} &
   &
  \multicolumn{1}{l}{} &
  \multicolumn{1}{l}{} \\
$\sigma^U=5$ &
  $\Delta^A$ &
   &
  1.21 &
  2.47 &
   &
  4.95 &
  6.65 &
   &
  10.35 &
  15.81 &
   &
  25.16 &
  60.65 \\
 &
  $\Delta^S$ &
   &
  0.04 &
  0.03 &
   &
  0.08 &
  0.07 &
   &
  0.14 &
  0.09 &
   &
  0.12 &
  0.11 \\ \bottomrule[1.5pt]
\end{tabular}%
}
\end{table}

Table \ref{tab:compare-CAOP} presents the average values of $\Delta^A$ and $\Delta^S$ across different parameter settings and problem sizes. Here, `N.A.' means that \citet{aouad2023exponomial} algorithm fails to terminate within 1 hour.  Overall, \citet{aouad2023exponomial} algorithm slightly outperforms our SBBD when both the utility variance $\sigma^U$ and the problem size are small, i.e., when $\sigma^U \in \{1,2\}$ and $|\ml{ J}| = \{10,20\}$. In these cases, \citet{aouad2023exponomial} algorithm manages to find  optimal solutions for nearly all tested instances. Our SBBD also demonstrates strong performance, as indicated by the small values of $\Delta^S$. 

However, we observe significant increases in $\Delta^A$ when the utility variance $\sigma^U$ increases. In particular, under $(\sigma^U, \sigma^R,|\ml{J}|) = (5,0.5,40)$, $\Delta^A$ reaches 60.65\%, implying that the objective provided by \citet{aouad2023exponomial} algorithm only achieves around 40\% of the upper bound on average. On the contrary, $\Delta^S$ exhibits slight increases with the maximum $\Delta^S$ being only 0.39\%, indicating that our SBBD obtains almost optimal solutions.

Besides being more stable in finding high-quality solutions,  our SBBD can solve larger problems. For instance,  \citet{aouad2023exponomial} algorithm is unable to solve instances with $|\ml{J}| = 40$ within 1 hour.  Moreover, when we ran \citet{aouad2023exponomial} algorithm in instances with $|\mathcal{J}| \ge 45$, our device runs out of memory without producing any solutions.  

\JJ{To further assess the efficiency of SBBD, we conducted computational experiments for 24 instances with $|\mathcal{A}| = 200$ ($|\ml{J}| = 201$) products and $\gamma = 0.1$ and compared its performance with respect to MILP and BIBC. Note that we use a stabilization procedure as described in Section \ref{sec:implement}. We construct the initial $x^s$ such that $x^s_0 = 1, x^s_j = \tau / |\ml{A}|,\forall j\in \ml{A}$. For the step size, we set $\varrho^{0} = 0.5$ and update $\varrho^{t} = \min\{\varrho^{t-1} + 0.05, 1\}$ at each iteration.}

\begin{table}[htb]
\centering
\caption{Efficiency test on CAOP under exponomial choice model.}
\label{tab:efficiency-expon}
\resizebox{\textwidth}{!}{%
\begin{tabular}{ccclccccclcccclccccc}
\toprule[1.5pt]
\multirow{2}{*}{$N$} &
  \multirow{2}{*}{$\sigma^R$} &
  \multirow{2}{*}{$\sigma^U$} &
   &
  \multicolumn{5}{c}{SBBD} &
   &
  \multicolumn{4}{c}{MILP} &
   &
  \multicolumn{5}{c}{BIBC} \\ \cline{5-9} \cline{11-14} \cline{16-20} 
 &
   &
   &
   &
  t{[}s{]} &
  \#Node &
  \#C &
  rgap{[}\%{]} &
  ogap{[}\%{]} &
   &
  t{[}s{]} &
  \#Node &
  rgap{[}\%{]} &
  ogap{[}\%{]} &
   &
  t{[}s{]} &
  \#Node &
  \#C &
  rgap{[}\%{]} &
  ogap{[}\%{]} \\ \hline
100 & 0.2 & 0.5 &  & 47.8   & 2053  & 10805 & 3.40  & 0.00 &  & 97.1   & 406 & 3.39  & 0.00 &  & 212.2  & 19136  & 4640  & 6.13  & 0.00  \\
    &     & 1   &  & 48.5   & 3922  & 11358 & 3.39  & 0.00 &  & 86.3   & 506 & 6.47  & 0.00 &  & 211.9  & 15334  & 5541  & 6.13  & 0.00  \\
    &     & 1.5 &  & 86.5   & 3623  & 9129  & 6.93  & 0.00 &  & 81.9   & 566 & 6.47  & 0.00 &  & 535.0  & 41675  & 7372  & 10.90 & 0.00  \\
    &     & 2   &  & 108.3  & 10044 & 3340  & 10.05 & 0.00 &  & 84.4   & 688 & 9.46  & 0.00 &  & 1343.6 & 102340 & 8299  & 15.00 & 0.00  \\
 &
  0.5 &
  0.5 &
   &
  12.2 &
  1745 &
  6839 &
  5.05 &
  0.00 &
   &
  101.1 &
  263 &
  4.79 &
  0.00 &
  \multicolumn{1}{c}{} &
  175.60 &
  15089 &
  4542 &
  9.01 &
  0.00 \\
 &
   &
  1 &
   &
  26.3 &
  2115 &
  8074 &
  7.12 &
  0.00 &
   &
  84.9 &
  270 &
  5.81 &
  0.00 &
  \multicolumn{1}{c}{} &
  143.10 &
  9911 &
  6059 &
  9.64 &
  0.00 \\
    &     & 1.5 &  & 39.3   & 2535  & 8608  & 7.62  & 0.00 &  & 68.9   & 274 & 7.46  & 0.00 &  & 153.7  & 14789  & 6590  & 14.92 & 0.00  \\
    &     & 2   &  & 42.5   & 2318  & 6935  & 7.66  & 0.00 &  & 46.4   & 173 & 2.69  & 0.00 &  & 167.7  & 9671   & 6905  & 20.47 & 0.00  \\
300 & 0.2 & 0.5 &  & 343.0  & 3741  & 27706 & 5.74  & 0.00 &  & 1365.0 & 596 & 5.08  & 0.00 &  & 3600.0 & 35397  & 15292 & 7.51  & 1.36  \\
    &     & 1   &  & 654.4  & 3824  & 38780 & 5.74  & 0.00 &  & 1157.1 & 521 & 5.07  & 0.00 &  & 3054.8 & 53725  & 19227 & 9.21  & 0.00  \\
    &     & 1.5 &  & 685.5  & 5662  & 32487 & 5.53  & 0.00 &  & 611.2  & 339 & 4.08  & 0.00 &  & 3600.0 & 32252  & 22926 & 10.47 & 1.18  \\
    &     & 2   &  & 839.6  & 12842 & 61713 & 9.39  & 0.00 &  & 920.5  & 713 & 8.99  & 0.00 &  & 3600.0 & 58937  & 23769 & 16.48 & 2.21  \\
    & 0.5 & 0.5 &  & 105.9  & 831   & 40381 & 4.08  & 0.00 &  & 1143.5 & 271 & 3.94  & 0.00 &  & 1920.7 & 7211   & 12249 & 7.44  & 0.00  \\
    &     & 1   &  & 201.4  & 1955  & 24092 & 5.29  & 0.00 &  & 781.1  & 271 & 4.17  & 0.00 &  & 1124.6 & 9788   & 16240 & 11.98 & 0.00  \\
    &     & 1.5 &  & 243.2  & 2686  & 33203 & 7.09  & 0.00 &  & 716.6  & 505 & 6.93  & 0.00 &  & 1932.0 & 17636  & 20651 & 18.38 & 0.00  \\
    &     & 2   &  & 290.7  & 5012  & 27093 & 12.52 & 0.00 &  & 468.1  & 291 & 10.60 & 0.00 &  & 3600.0 & 26082  & 20907 & 25.40 & 4.19  \\
500 & 0.2 & 0.5 &  & 804.1  & 4585  & 43497 & 5.29  & 0.00 &  & 3600.0 & 157 & 3.88  & 2.52 &  & 3600.0 & 21425  & 30110 & 8.33  & 4.04  \\
    &     & 1   &  & 971.3  & 9361  & 55691 & 4.93  & 0.00 &  & 3600.0 & 478 & 4.88  & 1.93 &  & 3600.0 & 18478  & 42053 & 17.29 & 14.10 \\
    &     & 1.5 &  & 1403.8 & 6840  & 62285 & 7.53  & 0.00 &  & 3146.9 & 455 & 6.57  & 0.00 &  & 3600.0 & 19673  & 38179 & 13.26 & 5.03  \\
    &     & 2   &  & 1237.9 & 11948 & 66911 & 10.16 & 0.00 &  & 3600.0 & 619 & 9.67  & 3.32 &  & 3600.0 & 13588  & 42539 & 22.03 & 18.35 \\
    & 0.5 & 0.5 &  & 242.2  & 1792  & 31780 & 4.72  & 0.00 &  & 3600.0 & 150 & 4.33  & 3.99 &  & 3600.0 & 10748  & 21445 & 10.63 & 0.79  \\
    &     & 1   &  & 431.0  & 2188  & 33138 & 5.86  & 0.00 &  & 3600.0 & 289 & 5.13  & 3.63 &  & 3600.0 & 13198  & 26275 & 11.33 & 3.77  \\
    &     & 1.5 &  & 1393.8 & 5058  & 69386 & 9.43  & 0.00 &  & 3600.0 & 113 & 7.90  & 4.14 &  & 3600.0 & 19282  & 31815 & 20.41 & 4.52  \\
    &     & 2   &  & 1154.8 & 4468  & 64167 & 12.10 & 0.00 &  & 2178.9 & 282 & 9.57  & 0.00 &  & 3600.0 & 9171   & 32343 & 25.33 & 16.00 \\ \bottomrule[1.5pt]
\end{tabular}%
}
\end{table}

\JJ{Table \ref{tab:efficiency-expon} reports the results. SBBD successfully solves all 24 instances within 1h; MILP solves 18 instances while BIBC solves only 12 instances. SBBD stands out because it exhibits the lowest minimum, maximum and average computational times. On average, BIBC explores much more nodes than SBBD and MILP, probably because of its weak relaxation as indicated by the root gap. Although MILP is able to maintain a small number of explored nodes for all the tested instances, it is still inferior than SBBD. An interesting observation is that when $\sigma^U$ increases, the computational time and root gap tend to increase for both SBBD and BIBC.} 

To summarize, in terms of solution quality, our SBBD has more stable performance and significantly outperforms \citet{aouad2023exponomial} algorithm when the utility variance $\sigma^U$ is large. In terms of computational speed, our SBBD is more efficient than the benchmarks and is capable of solving large instances.

\subsubsection{CAOP under mixed multinomial logit model}

We then conduct experiments on CAOP under mixed multinomial logit model (MMNL). We adopt the data generation method described in \cite{rusmevichientong2014assortment} to create our testbed. Specifically, the utility of a product $a \in \ml{A}$ is given by $U_a(\bs{\xi}) = V_a(\tilde{\boldsymbol{\xi}}) +  \epsilon_a$, and for the outside option, we set $U_0 = \epsilon_{0}$. Here, $\bs{\xi} = (\bs{\tilde{\xi}}, \bs{\epsilon})$, where $\tilde{\boldsymbol{\xi}}$ reflects the uncertainty of the intrinsic (or ideal) utility that contains latent customer segmentation information \footnote{For example, we can write $V_a(\tilde{\boldsymbol{\xi}}) = \tilde{\boldsymbol{\xi}}^T\bs{z}_a$, where $\bs{z}_a$ is the feature vector. In this case, we can interpret $\tilde{\boldsymbol{\xi}}$ as random sensitivity vector, and customer segments can be sampled from the distribution of $\tilde{\boldsymbol{\xi}}$.}, and $\bs{\epsilon}$ shows customer' idiosyncratic preferences and is assumed to follow i.i.d. Gumbel distribution.  According to~\cite{rusmevichientong2014assortment}, the uncertainty in the ideal utility of a product is captured by two random parameters $\tilde{\boldsymbol{\xi}}_a = (\tilde{\beta}_a, \tilde{\gamma}_a), \forall a \in \ml{A}$, and $V_a$ is defined as
\begin{equation}
	V_a = \ln(\tilde{\alpha}_a \tilde{\beta} _a) + \vartheta_1 + \vartheta_2\tilde{\gamma}_a , \forall a \in \ml{A}
\end{equation}
where  $\tilde{\beta}_a$ is a Bernoulli random variable with parameter 1/2, and $\tilde{\gamma}_a$ is a standard Normal random variable. The other parameters in $V_a$ are deterministic and generated as follows: $\tilde{\alpha}_a$ is the inherent attractiveness of option $a\in \mathcal{A}$, derived by first sampling $\psi_a$ from the standard normal distribution and then setting $\tilde{\alpha}_a = \psi_a / \sum_{a\in \mathcal{A}}\psi_a$. $\vartheta_1$ and $\vartheta_2$ depend on $D$ such that $\vartheta_1 + \vartheta^2_2 / 2 = \log10$, $\sqrt{e^{\vartheta^2_2 / 2} - 1} = D$. Moreover, the reward $R_a$ is sampled from a uniform distribution on the interval [1, $\bar{R}$] and $R_0 = 0$ for the outside option~\footnote{We refer the interested reader to \cite{rusmevichientong2014assortment} for understanding the reasoning behind this specific instance generation process.}.  

Under mixed multinomial logit model (MMNL), CAOP can be approximated as
\begin{align}
\max\quad &\frac{1}{\tilde{N}}\sum_{i\in \mathcal{I}}\frac{\sum_{a\in\mathcal{A}}R_{a}x_ae^{V_{ia}}}{\sum_{a\in\mathcal{A}} x_a e^{V_{ia}} + 1}\label{obj:MMNL}\\
\textbf{[CAOP-MMNL]}\qquad \text{s.t.}\quad &\sum_{a\in\mathcal{A}}x_a \le \tau,\\
 & x_a\in\{0,1\},\forall a\in\mathcal{A}.
\end{align}
where $\mathcal{I} = \{1,2,..,\tilde{N}\}$ represents the set of scenarios generated by sampling $\tilde{\boldsymbol{\xi}}$. The expectation over $\boldsymbol{\epsilon}$ is explicitly formulated in the linear fractional form, which avoids the need for sampling $\boldsymbol{\epsilon}$. To avoid confusion, we use a different notation $\tilde{N}$ to denote the sample size for \textbf{[CAOP-MMNL]}

The current state-of-the-art algorithm for solving \textbf{[CAOP-MMNL]} is a conic integer programming approach with McCormick estimators (conic+MC) proposed by \cite{sen2018conic}. We use this algorithm, which we call conic+MC, as a benchmark to compare with our SBBD.

Note that our SBBD samples both $\tilde{\boldsymbol{\xi}}$  and $\boldsymbol{\epsilon}$ (using LHS). Thus, our SBBD involves higher dimension of uncertainties and may require larger sample sizes than conic+MC. 
To reflect the trade-off between the solution quality and solution speed, we generate a different number of scenarios for conic+MC and compare its performance with our SBBD.

We fix $\tau = 5$ and vary $|\mathcal{A}| \in \{10,20,30\}, \bar{R} \in \{10,100,1000\}, D\in\{2,4,6\}$ to generate 27 instances. We use $\boldsymbol{x^C}$ and $\boldsymbol{x^S}$ to denote the solutions found by conic+MC and SBBD, respectively. The true objective values $v^C$ and $v^S$ are then estimated by a larger sample with $10^6$ scenarios. We compare the solution quality using the index $\textit{Diff}[\%] = (v^C - v^S) / v^C \times 100\%$. A negative value of \textit{Diff} indicates that SBBD finds better solution than conic+MC, and vice versa.

\begin{table}[h]
\centering
\caption{Comparison of conic+MC and SBBD on CAOP under MMNL with different sample sizes. `N.A.' indicates that conic+MC fails to obtain a feasible solution within 1 hour.}
\label{tab:compare-MMNL}
\resizebox{\textwidth}{!}{%
\begin{tabular}{rrrlrrrlrrrrlrrrrlrrrr}
\toprule[1.5pt]
$|\mathcal{A}|$ &
  $\bar{R}$ &
  $D$ &
   &
  \multicolumn{3}{c}{SBBD($N = 500$)} &
   &
  \multicolumn{4}{c}{conic+MC($\tilde{N} = 50$)} &
   &
  \multicolumn{4}{c}{conic+MC($\tilde{N} = 100$)} &
   &
  \multicolumn{4}{c}{conic+MC($\tilde{N} = 500$)} \\ \cline{5-7} \cline{9-12} \cline{14-17} \cline{19-22} 
 &
   &
   &
  \multicolumn{1}{c}{} &
  t{[}s{]} &
  \#Node &
  $v^S$ &
   &
  t{[}s{]} &
  \#Node &
  $v^C$ &
\textit{Diff}[\%]  &
   &
  t{[}s{]} &
  \#Node &
  $v^C$ &
\textit{Diff}[\%]  &
   &
  t{[}s{]} &
  \#Node &
  $v^C$ &
\textit{Diff}[\%] \\ \hline
10 & 10   & 2 &  & 0.15 & 1  & 2.28   &  & 0.10  & 1 & 2.18   & -4.62  &  & 3.61   & 1 & 2.28   & 0.00  &  & 27.59   & 1    & 2.28   & 0.00 \\
   &      & 4 &  & 0.13 & 1  & 1.61   &  & 0.07 & 1 & 1.53   & -5.25  &  & 4.04   & 1 & 1.61   & 0.00  &  & 175.94  & 7    & 1.61   & 0.00 \\
   &      & 6 &  & 0.13 & 1  & 1.28   &  & 0.07 & 1 & 1.22   & -5.23  &  & 7.27   & 1 & 1.28   & 0.00  &  & 583.22  & 29   & 1.28   & 0.00 \\
   & 100  & 2 &  & 0.14 & 1  & 21.87  &  & 0.09 & 1 & 20.98  & -4.25  &  & 3.85   & 1 & 21.87  & 0.00  &  & 28.69   & 1    & 21.87  & 0.00 \\
   &      & 4 &  & 0.12 & 1  & 15.44  &  & 0.07 & 1 & 14.74  & -4.77  &  & 5.84   & 1 & 15.44  & 0.00  &  & 27.26   & 1    & 15.44  & 0.00 \\
   &      & 6 &  & 0.11 & 1  & 12.31  &  & 0.07 & 1 & 11.73  & -4.97  &  & 4.89   & 1 & 12.31  & 0.00  &  & 120.77  & 3    & 12.31  & 0.00 \\
   & 1000 & 2 &  & 0.14 & 1  & 217.78 &  & 0.08 & 1 & 209.07 & -4.17  &  & 4.29   & 1 & 217.78 & 0.00  &  & 32.48   & 1    & 217.78 & 0.00 \\
   &      & 4 &  & 0.13 & 1  & 153.77 &  & 0.07 & 1 & 146.84 & -4.72  &  & 4.03   & 1 & 153.77 & 0.00  &  & 26.02   & 1    & 153.77 & 0.00 \\
   &      & 6 &  & 0.12 & 1  & 122.60 &  & 0.07 & 1 & 116.86 & -4.92  &  & 6.33   & 1 & 122.60 & 0.00  &  & 36.46   & 1    & 122.60 & 0.00 \\
20 & 10   & 2 &  & 0.20 & 1  & 2.56   &  & 0.17 & 1 & 2.37   & -8.02  &  & 46.94  & 1 & 2.65   & 3.40  &  & 2312.89 & 1    & 2.65   & 3.40 \\
   &      & 4 &  & 0.15 & 1  & 1.92   &  & 0.14 & 1 & 1.66   & -15.83 &  & 30.22  & 1 & 1.92   & 0.00  &  & 1083.57 & 1    & 1.92   & 0.00 \\
   &      & 6 &  & 0.12 & 1  & 1.56   &  & 0.15 & 1 & 1.34   & -16.37 &  & 24.47  & 1 & 1.56   & 0.00  &  & 696.74  & 1    & 1.56   & 0.00 \\
   & 100  & 2 &  & 0.20 & 1  & 24.25  &  & 0.17 & 1 & 22.34  & -8.56  &  & 47.42  & 1 & 25.27  & 4.03  &  & 1713.37 & 1    & 25.27  & 4.03 \\
   &      & 4 &  & 0.11 & 1  & 18.37  &  & 0.22 & 1 & 15.96  & -15.12 &  & 23.9   & 1 & 18.37  & 0.00  &  & 727.90  & 1    & 18.37  & 0.00 \\
   &      & 6 &  & 0.11 & 1  & 14.91  &  & 0.15 & 1 & 11.79  & -26.48 &  & 30.62  & 1 & 14.91  & 0.00  &  & 642.78  & 1    & 14.91  & 0.00 \\
   & 1000 & 2 &  & 0.16 & 1  & 244.17 &  & 0.16 & 1 & 222.05 & -9.96  &  & 46.23  & 1 & 251.45 & 2.89  &  & 1801.48 & 1    & 251.45 & 2.89 \\
   &      & 4 &  & 0.11 & 1  & 182.87 &  & 0.14 & 1 & 158.62 & -15.29 &  & 125.44 & 1 & 182.87 & 0.00  &  & 736.71  & 1    & 182.87 & 0.00 \\
   &      & 6 &  & 0.11 & 1  & 148.43 &  & 0.14 & 1 & 117.15 & -26.70 &  & 27.93  & 1 & 148.43 & 0.00  &  & 646.22  & 1    & 148.43 & 0.00 \\
30 & 10   & 2 &  & 0.23 & 11 & 2.81   &  & 0.26 & 1 & 2.62   & -7.37  &  & 615.24 & 1 & 2.81   & 0.00  &  & 3600.00 & N.A. & N.A.   & N.A. \\
   &      & 4 &  & 0.20 & 3  & 2.01   &  & 0.14 & 1 & 1.9    & -5.92  &  & 373.93 & 1 & 1.93   & -4.27 &  & 3600.00 & N.A. & N.A.   & N.A. \\
   &      & 6 &  & 0.18 & 1  & 1.66   &  & 0.15 & 1 & 1.54   & -7.48  &  & 63.97  & 1 & 1.56   & -6.10 &  & 3389.82 & 1    & 1.66   & 0.00 \\
   & 100  & 2 &  & 0.22 & 5  & 26.69  &  & 0.29 & 1 & 27.01  & 1.20   &  & 609.38 & 3 & 27.01  & 1.20  &  & 3600.00 & N.A. & N.A.   & N.A. \\
   &      & 4 &  & 0.16 & 1  & 19.37  &  & 0.16 & 1 & 17.99  & -7.69  &  & 336.22 & 1 & 18.47  & -4.89 &  & 3600.00 & N.A. & N.A.   & N.A. \\
   &      & 6 &  & 0.16 & 1  & 15.91  &  & 0.16 & 1 & 14.55  & -9.33  &  & 122.81 & 1 & 14.96  & -6.33 &  & 3600.00 & N.A. & N.A.   & N.A. \\
   & 1000 & 2 &  & 0.26 & 3  & 266.07 &  & 0.28 & 1 & 269.09 & 1.12   &  & 673.01 & 3 & 269.09 & 1.12  &  & 3600.00 & N.A. & N.A.   & N.A. \\
   &      & 4 &  & 0.18 & 1  & 192.99 &  & 0.15 & 1 & 178.82 & -7.92  &  & 384.94 & 1 & 183.89 & -4.95 &  & 3600.00 & N.A. & N.A.   & N.A. \\
   &      & 6 &  & 0.18 & 3  & 156.48 &  & 0.16 & 1 & 144.67 & -8.16  &  & 153.97 & 1 & 148.88 & -5.10 &  & 3600.00 & N.A. & N.A.   & N.A. \\ \bottomrule[1.5pt]
\end{tabular}%
}
\end{table}

Table \ref{tab:compare-MMNL} reports the result. Here, we generate $\tilde{N}\in\{50,100,500\}$ scenarios for conic+MC and $N = 500$ scenarios for SBBD. `N.A.'  indicates that the approach fails to obtain a feasible solution within the 1-hour time limit. We find that our SBBD is able to solve all the tested instances within 0.3s. Under $\tilde{N} = 50$ scenarios, the speed of conic+MC is comparable with SBBD. However, conic+MC exhibits poor solution quality. On average, the revenue by conic+MC($\tilde{N} = 50$) is 8.77\% worse than that of SBBD, and the worst case is 26.70\%. When $\tilde{N}$ increases to 100, the solution quality of conic+MC greatly improves. The average difference is -0.69\%, meaning that conic+MC has comparable solution quality with SBBD. However, it comes at the cost of increased computational burden, as the average computational time of conic+MC becomes 140s, much higher than that of SBBD. Finally, when $\tilde{N} = 500$, conic+MC can provide the same or even better solutions than SBBD when $|\mathcal{A}|\in\{10,20\}$. However, when $|\mathcal{A}| = 30$, conic+MC is only able to solve one instance within 1h. Moreover, for the remaining 8 instances, conic+MC cannot even find a feasible solution. 

In summary, our SBBD with $N = 500$ scenarios can efficiently find high quality solutions for CAOP under MMNL. With $N = 100$ scenarios, the solution quality of conic+MC is comparable with SBBD, because of not requiring sampling $\boldsymbol{\epsilon}$. However, the average computational time of conic+MC is significantly longer than SBBD. When we can further increase $N$ to 500 for conic+MC, the computational burden makes it very inefficient for solving the problem. This underscores SBBD's superior capabilities in balancing between solution quality and computational time.

\remark In addition to the exponomial choice model and mixed logit model, we also conduct experiments with the multinomial probit model, which lacks a closed-form choice expression. The results affirm that our methodology yields strong performance for CAOP. Due to page constraints, we have relegated the results to Online Appendices A.2 and A.3.

\subsection{Experiment on FLoP and MSMFLP}
In this section, we handle two facility location problems. We generate synthetic data as follows.

\textbf{FLoP}. The geographical distribution of customers is assumed to follow a uniform distribution with $\bs{\xi} = (lon,lat) \sim\mathcal{U}[0, 20] \times \mathcal{U}[0,20]$.  Therefore, we sample a set of customers, $\mathcal{I}$, from the distribution. For the candidate facilities $a\in\mathcal{A}$, we generate their coordinates  from $\mathcal{U}[0,20] \times \mathcal{U}[0,20]$. To calculate the distance $d_{ia}$ between customer at $i\in\mathcal{I}$ and a facility $a\in\mathcal{A}$, we employ the Euclidean distance formula. For the sake of our analysis, we consider 10 pricing levels (i.e., $|\mathcal{L}| = 10$), where each level $l$ is associated with a pricing value $s_l = l$ for $l = 1, 2, 3,..., 10$. Consequently, we define $G_{ial} = d_{ia} + s_l,\forall i\in\mathcal{I},\forall a\in\mathcal{A},\forall l\in\mathcal{L}$ (\JJ{i.e., $\bs{d}$ and $\bs{s}$ have the same marginal contribution to $\bs{G}$}). For the outside option, we fix $b = 10$.

\textbf{MSMFLP}. The geographical distribution of customers is assumed to follow a normal distribution with $\bs{\xi} = (lon,lat) \sim\mathcal{N}[10,100/3] \times \mathcal{N}[10,100/3]$. Therefore, we sample a set of customers, $\mathcal{I}$, from the distribution.
For the candidate facilities, we generate their coordinates from $\mathcal{U}[0,20] \times \mathcal{U}[0,20]$. To calculate the distance $d_{ij}$ between sampled customer $i\in\mathcal{I}$ and facility $j\in\mathcal{J}$, we employ the Euclidean distance formula. The attractions of facilities are generated as $a_j\sim\mathcal{U}[1,20],\forall j\in\mathcal{J}$. The utility of facility $j\in\mathcal{J}$ perceived by sampled customer $i\in\mathcal{I}$ is computed as $u_{ij} = a_j / d^2_{ij}$. The utility of the outside option is defined as $O = 10$.

\subsubsection{Solution quality test}

Since we are not aware of any other solution approaches that can directly solve the stochastic versions of FLoP and MSMFLP, we are unable to compare the solution quality as we did for CAOP.  Instead, we leverage the approach described in Section~\ref{sec:SAA} and use the estimated gaps $\Delta$ and $\Delta_{0.95}$ to measure the solution quality. Moreover, there are no closed-form formulations of the objective functions for FLoP and MSMFLP. Therefore, we use a large sample $N' = 1\times10^6$ to estimate $\hat{v}$, i.e., the true value of the objective obtained by our SBBD.

To generate the testbed, we set $|\mathcal{A}|\in\{10,20\}$ for FLoP, leading to $|\ml{J}| = |\ml{A}||\ml{L}| + 1\in\{101,201\}$. For MSMFLP, we use $|\mathcal{J}|\in\{100,200\}$ . We fix $\tau = 5$ for FLoP and $\tau = 20$ for MSMFLP. We vary sample size $N\in \{300,500,1000\}$ and replicate $M = 50$ times by changing the random seed from 1 to 50 when sampling $\bs{\xi}$. 

\begin{figure}[htb]
	\centering
	\subfigure[$\Delta$ of FLoP]{
		\psfig{figure=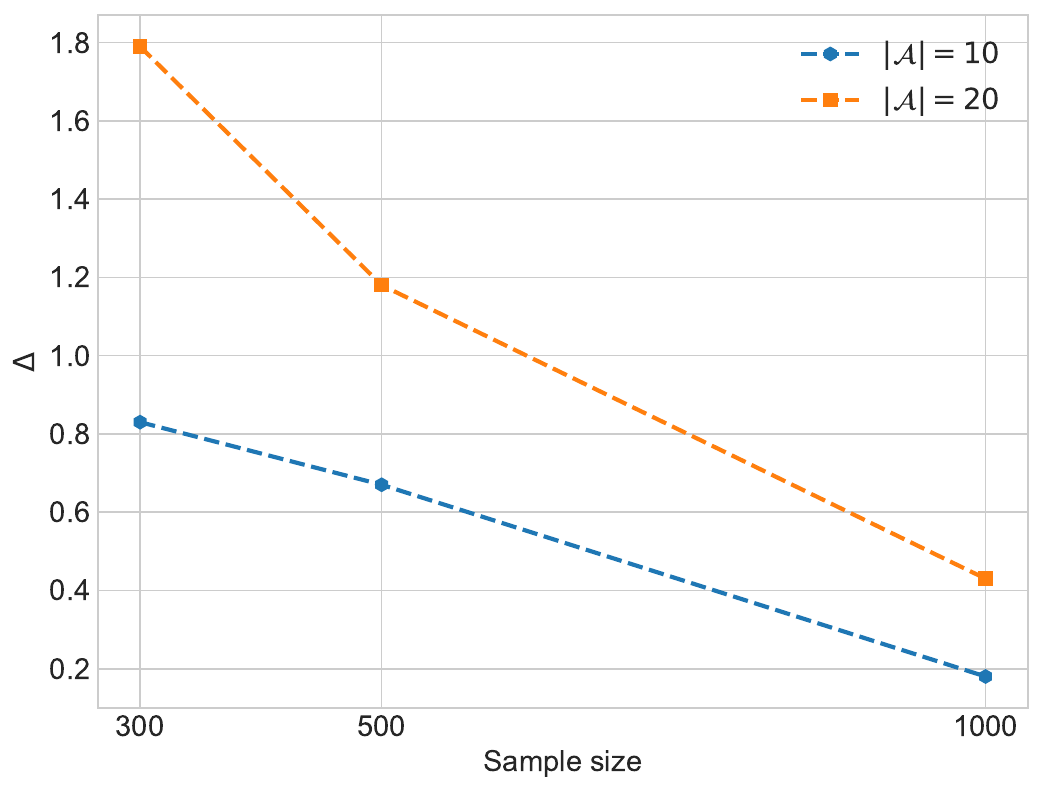, width=77mm, height=55mm}
		\label{fig:app2-D2}
	}
	\subfigure[$\Delta_{0.95}$ of FLoP]{
		\psfig{figure=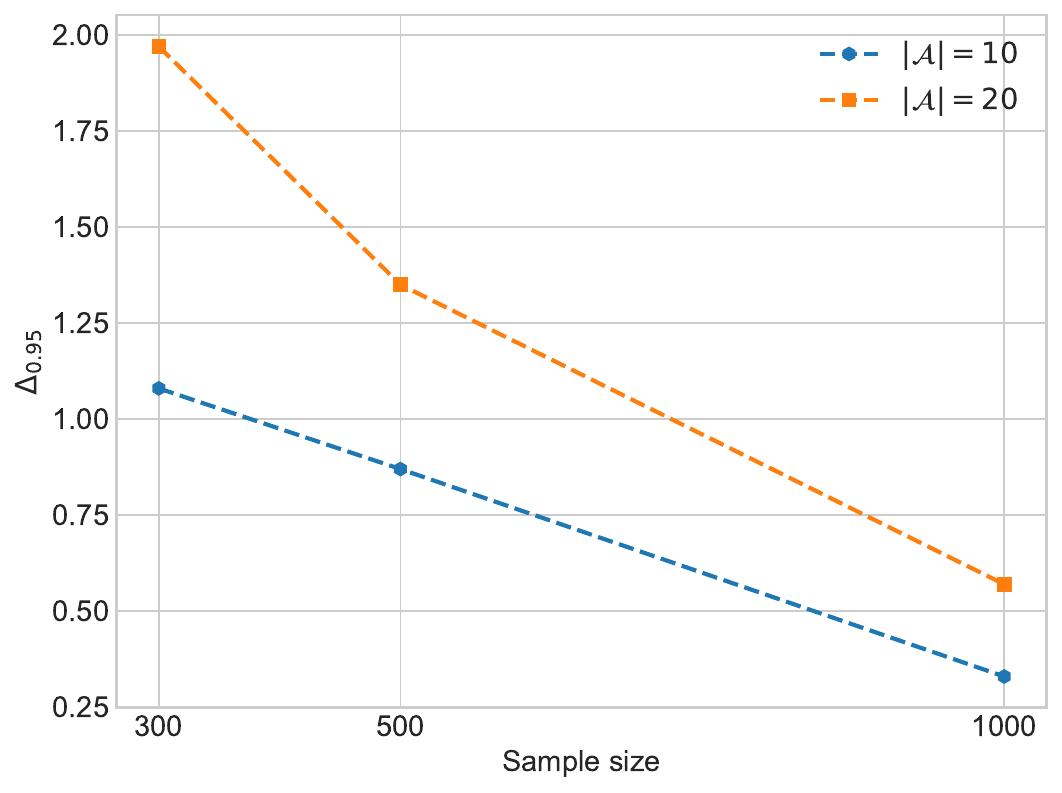, width=77mm, height=55mm}
		\label{fig:app2-Da2}
	}

	\subfigure[$\Delta$ of MSMFLP]{
		\psfig{figure=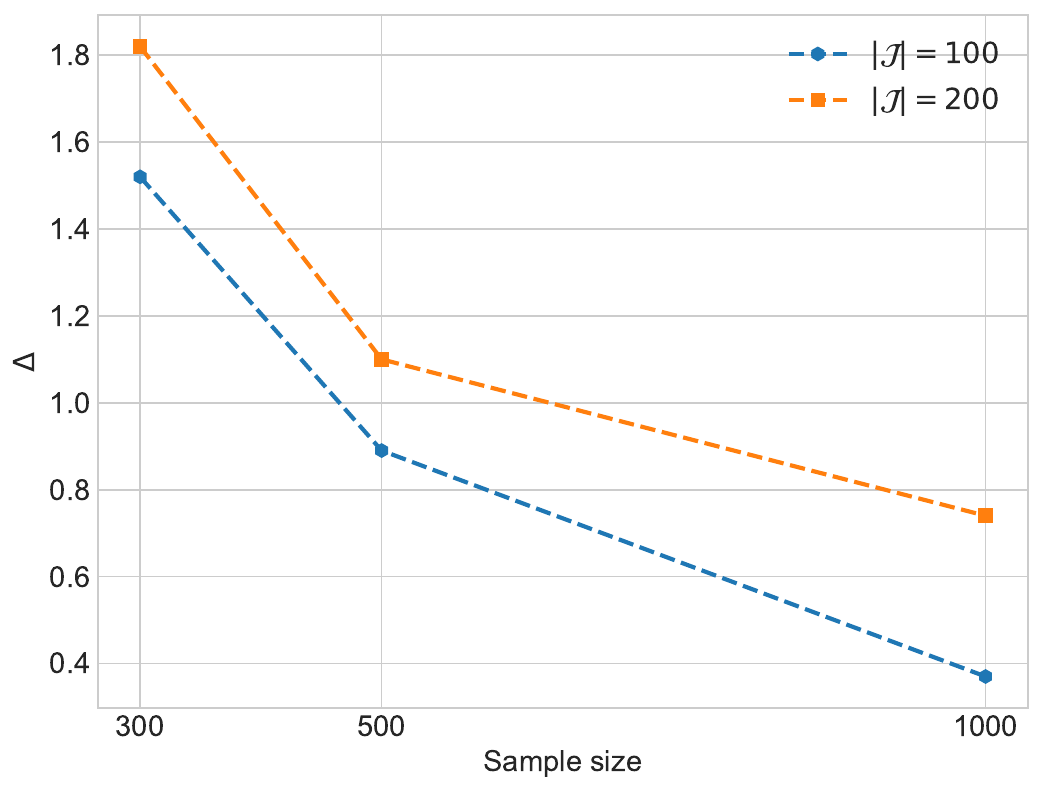, width=77mm, height=55mm}
		\label{fig:app3-D3}
	}
	\subfigure[$\Delta_{0.95}$ of MSMFLP]{
		\psfig{figure=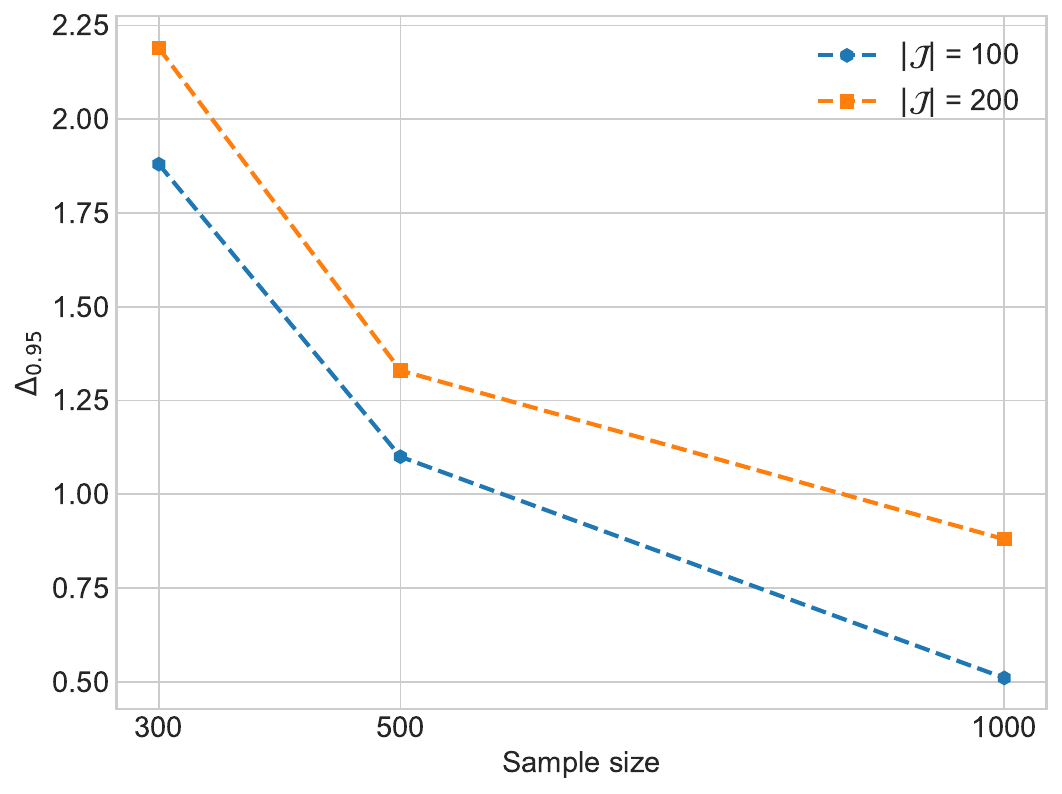, width=77mm, height=55mm}
		\label{fig:app3-Da3}
	}
	\caption{Estimated gaps $\Delta$ and $\Delta_{0.95}$ for FLoP and MSMFLP}
	\label{fig:sol-val}
\end{figure}

The values of $\Delta$ and $\Delta_{0.95}$ are illustrated in Figure \ref{fig:sol-val}. The complete results  are available in Online Appendix A.2.2. Intuitively, both $\Delta$ and $\Delta_{0.95}$ decrease as the sample size increases. Additionally, we observe that larger-scale (larger $|\ml{J}|$) instances typically show higher values for $\Delta$ and $\Delta_{0.95}$. However, we only observe very limited increase of $\Delta$ and $\Delta_{0.95}$ regardless of sample sizes, consistently staying well below 1\% for both FLoP and MSMFLP (see Tables A.2 and A.3 in Online Appendix). 
When $N = 1000$, SBBD achieves exceptionally high solution quality for both problems, with both $\Delta$ and $\Delta_{0.95}$ dropping below 1\%. Indeed, a sample size of $N = 500$ proves to be sufficiently large, leading to $\Delta$ below 1\% for small-scale instances and approximately 1\% for larger-scale instances. This demonstrates that our SBBD is capable of generating high-quality solutions with reasonable sample sizes.

\subsubsection{Efficiency test}
We then proceed to test the efficiency of SBBD, using instances of larger scales than those used in the solution quality test.

\noindent \textbf{FLoP}. We generate 21 instances with varying scales and $\tau$. For benchmark solution approaches, we utilize MILP and BIBC. Table \ref{tab:efficiency-FLoP} presents the computational results for these instances. We observe that all the three approaches optimally solve all instances within 1h. Notably, SBBD demonstrates the best performance: it exhibits the lowest average computational time, which is just 26.43s, significantly shorter than that of MILP (239.16s) and BIBC (132.14s).

\begin{table}[htb]
\centering
\caption{Computational results of 21 FLoP instances}
\label{tab:efficiency-FLoP}
\resizebox{\textwidth}{!}{%
\begin{tabular}{cccccccccccccccccccc}
\toprule[1.5pt]
\multirow{2}{*}{$|\mathcal{A}|$} &
  \multirow{2}{*}{$\tau$} &
  \multirow{2}{*}{$N$} &
   &
  \multicolumn{5}{c}{SBBD} &
   &
  \multicolumn{4}{c}{MILP} &
   &
  \multicolumn{5}{c}{BIBC} \\ \cline{5-9} \cline{11-14} \cline{16-20} 
 &
   &
   &
   &
  t{[}s{]} &
  \#Node &
  \#C &
  rgap{[}\%{]} &
  ogap{[}\%{]} &
   &
  t{[}s{]} &
  \#Node &
  rgap{[}\%{]} &
  ogap{[}\%{]} &
   &
  t{[}s{]} &
  \#Node &
  \#C &
  rgap{[}\%{]} &
  ogap{[}\%{]} \\ \hline
20 & 5  & 300  &  & 2.4   & 14   & 2030  & 1.22 & 0.00 &  & 9.77    & 1   & 0.00 & 0.00 &  & 4.3    & 7    & 253  & 1.46 & 0.00 \\
   &    & 500  &  & 3.6   & 10   & 3602  & 0.98 & 0.00 &  & 22.25   & 1   & 0.00 & 0.00 &  & 10.2   & 14   & 661  & 0.96 & 0.00 \\
   &    & 1000 &  & 9.0   & 54   & 8461  & 2.06 & 0.00 &  & 91.38   & 1   & 0.00 & 0.00 &  & 29.2   & 119  & 1315 & 1.95 & 0.00 \\
30 & 5  & 300  &  & 5.0   & 13   & 2909  & 0.84 & 0.00 &  & 38.72   & 1   & 0.00 & 0.00 &  & 8.7    & 17   & 389  & 0.72 & 0.00 \\
   &    & 500  &  & 7.2   & 12   & 3299  & 0.34 & 0.00 &  & 65.78   & 1   & 0.00 & 0.00 &  & 13.5   & 8    & 538  & 0.43 & 0.00 \\
   &    & 1000 &  & 11.3  & 3    & 6324  & 0.43 & 0.00 &  & 135.43  & 1   & 0.00 & 0.00 &  & 33.7   & 4    & 1441 & 0.41 & 0.00 \\
30 & 10 & 300  &  & 4.5   & 107  & 2752  & 1.09 & 0.00 &  & 39.91   & 1   & 0.00 & 0.00 &  & 16.6   & 73   & 700  & 1.07 & 0.00 \\
   &    & 500  &  & 9.7   & 487  & 4934  & 2.28 & 0.00 &  & 148.70  & 33  & 0.56 & 0.00 &  & 83.8   & 1123 & 2331 & 2.36 & 0.00 \\
   &    & 1000 &  & 42.4  & 933  & 16199 & 2.44 & 0.00 &  & 397.50  & 180 & 1.31 & 0.00 &  & 486.8  & 1790 & 4876 & 2.31 & 0.00 \\
40 & 5  & 300  &  & 5.9   & 93   & 2146  & 1.98 & 0.00 &  & 70.76   & 1   & 0.00 & 0.00 &  & 9.4    & 10   & 326  & 0.45 & 0.00 \\
   &    & 500  &  & 12.8  & 808  & 5630  & 2.98 & 0.00 &  & 140.00  & 1   & 0.00 & 0.00 &  & 19.6   & 13   & 632  & 0.62 & 0.00 \\
   &    & 1000 &  & 16.4  & 399  & 4702  & 2.73 & 0.00 &  & 350.92  & 1   & 0.00 & 0.00 &  & 45.6   & 36   & 1352 & 0.53 & 0.00 \\
40 & 10 & 300  &  & 17.8  & 2069 & 7536  & 2.89 & 0.00 &  & 107.90  & 49  & 0.18 & 0.00 &  & 39.4   & 315  & 1003 & 1.51 & 0.00 \\
   &    & 500  &  & 17.3  & 607  & 4740  & 2.13 & 0.00 &  & 287.64  & 52  & 0.82 & 0.00 &  & 135.4  & 1470 & 2524 & 1.95 & 0.00 \\
   &    & 1000 &  & 164.8 & 4729 & 38799 & 2.92 & 0.00 &  & 1022.45 & 128 & 1.65 & 0.00 &  & 1051.1 & 2635 & 6486 & 2.38 & 0.00 \\
50 & 5  & 300  &  & 9.8   & 55   & 2558  & 1.18 & 0.00 &  & 51.06   & 1   & 0.00 & 0.00 &  & 15.9   & 46   & 426  & 0.95 & 0.00 \\
   &    & 500  &  & 16.4  & 85   & 3754  & 1.37 & 0.00 &  & 181.82  & 1   & 0.00 & 0.00 &  & 34.5   & 48   & 737  & 1.33 & 0.00 \\
   &    & 1000 &  & 29.6  & 64   & 8174  & 1.03 & 0.92 &  & 509.33  & 17  & 0.17 & 0.00 &  & 88.4   & 55   & 1383 & 1.21 & 0.00 \\
50 & 10 & 300  &  & 9.8   & 1158 & 2736  & 2.06 & 0.00 &  & 132.20  & 44  & 0.53 & 0.00 &  & 29.5   & 165  & 574  & 1.40 & 0.00 \\
   &    & 500  &  & 27.0  & 2019 & 6275  & 1.99 & 0.00 &  & 199.37  & 1   & 0.00 & 0.00 &  & 60.8   & 106  & 1128 & 0.83 & 0.00 \\
   &    & 1000 &  & 132.2 & 5070 & 22226 & 2.87 & 0.00 &  & 1019.55 & 99  & 1.42 & 0.00 &  & 558.7  & 691  & 3199 & 1.91 & 0.00 \\ \bottomrule[1.5pt]
\end{tabular}%
}
\end{table}

\eat{
\clearpage

\begin{table}[htb]
\centering
\caption{Summary of statistics of measure for FLoP.}
\label{tab:sum2}
\resizebox{0.8\textwidth}{!}{%
\begin{tabular}{ccccccccccccc}
\toprule[1.5pt]
 & \multirow{7}{*}{} & \multicolumn{3}{c}{SBBD} & \multirow{7}{*}{} & \multicolumn{3}{c}{MILP} & \multirow{7}{*}{} & \multicolumn{3}{c}{BIBC} \\ \cline{3-5} \cline{7-9} \cline{11-13}
             &  & Min  & Max    & Avg     &  & Min   & Max     & Avg     &  & Min  & Max     & Avg     \\ \cline{1-1} \cline{3-5} \cline{7-9} \cline{11-13}
t{[}s{]}     &  & 2.42 & 164.82 & 26.43   &  & 48.13 & 3600.00 & 1257.78 &  & 4.34 & 1051.11 & 132.14  \\
\#Node          &  & 3    & 5070   & 894.71  &  & 1     & 589     & 97.91   &  & 4    & 2635    & 416.43  \\
\#C          &  & 2987 & 42922  & 9838.52 &  & 1224  & 26504   & 5614.09 &  & 253  & 6486    & 1536.86 \\
rgap{[}\%{]} &  & 0.34 & 2.98   & 1.80    &  & 0.02  & 4.94    & 1.24    &  & 0.41 & 2.38    & 1.27    \\
ogap{[}\%{]} &  & 0.00 & 0.00   & 0.00    &  & 0.00  & 2.10    & 0.15    &  & 0.00 & 0.00    & 0.00    \\ \bottomrule[1.5pt]
\end{tabular}%
}
\end{table}

\begin{figure}[htb]
		\begin{center}
			\subfigure[Number of instances solved versus computational time]{
				\psfig{figure = figures/app2-cum, width=77mm, height = 55mm}\label{fig:app2-cum}}
			\subfigure[Boxplot of logarithm of computational time]{
				\psfig{figure = figures/app2-bp, width=77mm, height = 55mm}\label{fig:app2-bp} }
		\end{center}
		\caption{Time comparison of different solution approaches for FLoP.}
		\label{fig:app2}
\end{figure}
}

\noindent \textbf{MSMFLP}. The next experiment focuses on MSMFLP, for which the resulting model after sampling can be reformulated as a mixed-integer linear fractional program and can be solved by several solution approaches. In a recent paper by \citet{lin2021branch}, the authors propose a generalized Benders decomposition (GBD) approach for MSMFLP and show that it outperforms state-of-the-art approaches. Therefore, we use MILP and GBD as benchmarks for MSMFLP. 

\JJ{Table \ref{tab:efficiency-MSMFLP} shows the computational results of 30 MSMFLP instances of varying sizes and parameters. For the 18 instances with $|\ml{J}| = 200$, SBBD and MILP solves all of them within 1h, while GBD fails to solve one instance with optimality gap 0.9\%. This is because the weak relaxation of GBD, as indicated by the large root gaps as well as a remarkably large number of nodes  explored. For the remaining 12 instances with $|\ml{J}| = 400$, MILP runs out of memory and no solution is returned. Although GBD does not have a memory issue, it fails to solve all the instances with $|\ml{J}| = 400$. On the contrary, SBBD still solves all of them within 30s. }

\begin{table}[htb]
\centering
\caption{Computational results of 30 MSMFLP instances}
\label{tab:efficiency-MSMFLP}
\resizebox{\textwidth}{!}{%
\begin{tabular}{ccccccccccccccccccccc}
\toprule[1.5pt]
\multirow{2}{*}{$N$} &
  \multirow{2}{*}{$|\mathcal{J}|$} &
  \multirow{2}{*}{$p$} &
  \multirow{2}{*}{$O$} &
   &
  \multicolumn{5}{c}{SBBD} &
   &
  \multicolumn{4}{c}{MILP} &
   &
  \multicolumn{5}{c}{GBD} \\ \cline{6-10} \cline{12-15} \cline{17-21} 
 &
   &
   &
   &
   &
  t{[}s{]} &
  \#Node &
  \#C &
  rgap{[}\%{]} &
  ogap{[}\%{]} &
   &
  t{[}s{]} &
  \#Node &
  rgap{[}\%{]} &
  ogap{[}\%{]} &
   &
  t{[}s{]} &
  \#Node &
  \#C &
  rgap{[}\%{]} &
  ogap{[}\%{]} \\ \hline
500  & 200 & 10 & 3  &  & 2.7  & 1  & 4760  & 0.00 & 0.00 &  & 169.2 & 1    & 0.00 & 0.00 &  & 90.1   & 25671   & 9941   & 17.60 & 0.00 \\
     &     &    & 5  &  & 2.8  & 3  & 6310  & 0.05 & 0.00 &  & 257.8 & 1    & 0.00 & 0.00 &  & 112.0  & 67280   & 6781   & 12.88 & 0.00 \\
     &     &    & 7  &  & 2.6  & 1  & 4870  & 0.00 & 0.00 &  & 211.8 & 1    & 0.00 & 0.00 &  & 63.3   & 32952   & 7392   & 10.54 & 0.00 \\
     &     &    & 10 &  & 2.3  & 1  & 4570  & 0.00 & 0.00 &  & 393.9 & 1    & 0.00 & 0.00 &  & 45.1   & 24294   & 6771   & 8.09  & 0.00 \\
     &     &    & 15 &  & 2.4  & 3  & 5458  & 0.08 & 0.00 &  & 583.2 & 1    & 0.00 & 0.00 &  & 35.8   & 13286   & 7107   & 5.88  & 0.00 \\
     &     &    & 20 &  & 2.3  & 1  & 4483  & 0.00 & 0.00 &  & 323.8 & 1    & 0.00 & 0.00 &  & 40.6   & 13795   & 7926   & 4.20  & 0.00 \\
500  & 200 & 20 & 3  &  & 2.3  & 1  & 4496  & 0.00 & 0.00 &  & 121.5 & 1    & 0.00 & 0.00 &  & 769.8  & 1393181 & 6853   & 24.13 & 0.00 \\
     &     &    & 5  &  & 2.3  & 1  & 4425  & 0.00 & 0.00 &  & 118.8 & 1    & 0.00 & 0.00 &  & 1032.1 & 1772801 & 7318   & 17.65 & 0.00 \\
     &     &    & 7  &  & 2.2  & 1  & 4460  & 0.00 & 0.00 &  & 97.7  & 1    & 0.00 & 0.00 &  & 410.1  & 682577  & 8288   & 14.64 & 0.00 \\
     &     &    & 10 &  & 2.2  & 1  & 4432  & 0.00 & 0.00 &  & 103.0 & 1    & 0.00 & 0.00 &  & 332.9  & 705488  & 7159   & 11.17 & 0.00 \\
     &     &    & 15 &  & 2.2  & 1  & 4464  & 0.00 & 0.00 &  & 125.8 & 1    & 0.00 & 0.00 &  & 434.9  & 1071337 & 6840   & 8.11  & 0.00 \\
     &     &    & 20 &  & 2.7  & 4  & 5637  & 0.03 & 0.00 &  & 143.0 & 1    & 0.00 & 0.00 &  & 188.8  & 449683  & 6712   & 5.70  & 0.00 \\
1000 & 200 & 20 & 3  &  & 5.7  & 5  & 11616 & 0.03 & 0.00 &  & 306.8 & 1    & 0.00 & 0.00 &  & 3600.0 & 831934  & 17085  & 23.71 & 0.90 \\
     &     &    & 5  &  & 5.6  & 7  & 11092 & 0.02 & 0.00 &  & 378.0 & 1    & 0.00 & 0.00 &  & 1135.5 & 500820  & 15663  & 17.20 & 0.00 \\
     &     &    & 7  &  & 5.8  & 9  & 11644 & 0.01 & 0.00 &  & 423.3 & 1    & 0.00 & 0.00 &  & 563.1  & 221434  & 19459  & 14.96 & 0.00 \\
     &     &    & 10 &  & 4.5  & 1  & 16199 & 0.00 & 0.00 &  & 440.2 & 1    & 0.00 & 0.00 &  & 485.6  & 230668  & 16471  & 10.83 & 0.00 \\
     &     &    & 15 &  & 4.9  & 1  & 10573 & 0.00 & 0.00 &  & 411.5 & 1    & 0.00 & 0.00 &  & 386.2  & 213040  & 16726  & 8.68  & 0.00 \\
     &     &    & 20 &  & 4.5  & 1  & 8795  & 0.00 & 0.00 &  & 481.2 & 1    & 0.00 & 0.00 &  & 494.0  & 367442  & 13992  & 6.87  & 0.00 \\
1000 & 400 & 20 & 3  &  & 14.0 & 1  & 10884 & 0.00 & 0.00 &  & N.A.  & N.A. & N.A. & N.A. &  & 3600.0 & 176256  & 22700  & 29.62 & 4.34 \\
     &     &    & 5  &  & 12.5 & 6  & 11420 & 0.06 & 0.00 &  & N.A.  & N.A. & N.A. & N.A. &  & 3600.0 & 102709  & 33769  & 24.74 & 3.14 \\
     &     &    & 7  &  & 11.8 & 1  & 9576  & 0.00 & 0.00 &  & N.A.  & N.A. & N.A. & N.A. &  & 3600.0 & 253365  & 26071  & 20.02 & 2.31 \\
     &     &    & 10 &  & 12.0 & 1  & 9379  & 0.00 & 0.00 &  & N.A.  & N.A. & N.A. & N.A. &  & 3600.0 & 277096  & 26011  & 16.67 & 1.75 \\
     &     &    & 15 &  & 14.3 & 11 & 12180 & 0.03 & 0.00 &  & N.A.  & N.A. & N.A. & N.A. &  & 3600.0 & 242159  & 21766  & 7.09  & 1.31 \\
     &     &    & 20 &  & 11.6 & 1  & 9459  & 0.00 & 0.00 &  & N.A.  & N.A. & N.A. & N.A. &  & 3600.0 & 226858  & 231425 & 5.25  & 1.06 \\
2000 & 400 & 20 & 3  &  & 28.8 & 12 & 25150 & 0.08 & 0.00 &  & N.A.  & N.A. & N.A. & N.A. &  & 3600.0 & 21336   & 77812  & 31.84 & 8.22 \\
     &     &    & 5  &  & 26.7 & 34 & 27192 & 0.08 & 0.00 &  & N.A.  & N.A. & N.A. & N.A. &  & 3600.0 & 17521   & 66345  & 26.19 & 4.85 \\
     &     &    & 7  &  & 29.2 & 13 & 26616 & 0.06 & 0.00 &  & N.A.  & N.A. & N.A. & N.A. &  & 3600.0 & 26033   & 52788  & 19.63 & 3.76 \\
     &     &    & 10 &  & 28.2 & 29 & 24226 & 0.07 & 0.00 &  & N.A.  & N.A. & N.A. & N.A. &  & 3600.0 & 27991   & 117234 & 18.72 & 4.89 \\
     &     &    & 15 &  & 28.9 & 31 & 26289 & 0.09 & 0.00 &  & N.A.  & N.A. & N.A. & N.A. &  & 3600.0 & 34142   & 61006  & 16.43 & 2.52 \\
     &     &    & 20 &  & 25.9 & 30 & 25823 & 0.03 & 0.00 &  & N.A.  & N.A. & N.A. & N.A. &  & 3600.0 & 37596   & 74020  & 11.36 & 1.85 \\ \bottomrule[1.5pt]
\end{tabular}%
}
\end{table}

To summarize, our SBBD can obtain high-quality solutions for FLoP and MSMFLP with a reasonable number of sample size. Moreover, it significantly outperforms advanced solution approaches in the literature in terms of computational speed.

\section{Impact of the reward-utility relation}\label{sec:insight}
This section is inspired by the observation that MILP solves all the MSMFLP nstances at root nodes, and SBBD consistently achieves small root gaps (rgap <$ 0.1\%$) for all MSMFLP instances. This implies that our \textbf{[ADPP]} has a high level of formulation tightness. Remarkably, SBBD exhibits exceptional speed on MSMFLP, solving instances with up to  $|\ml{J}| = 400$ facilities and $N = 2000$ samples within 30s.

After sampling, \textbf{[ADPP]} provides a uniform formulation akin to the rank list choice model, where the customers have a fixed preference ranking over all options  -- an option's higher sampled utility translates to a higher preference ranking, leading customers to select the top-ranked option available to them. In this setting, the planner's objective is to maximize rewards, while customers seek to maximize their utilities. Notably, for MSMFLP, the ranking derived from sampled utilities aligns perfectly with the ranking based on sampled rewards, as highlighted in (\ref{R_on_xi}). This alignment implies that maximizing the utility of the customer leads to the maximization of the planner's reward, thereby simplifying the model's computational complexity.  Therefore, it is worth discussing the impact the reward-utility relation on computational difficulty of \textbf{[ADPP]}.

Consider a planning decision $\bs{x}\in\Omega$, for scenario $i\in\mathcal{I}$, let $j_i^r$ and $j_i^u$ be such that
\begin{align}
j_i^r = \mathop{\arg\max}_{j\in \ml{C}(\bs{x})}r_{ij}, \quad j_i^u = \mathop{\arg\max}_{j\in \ml{C}(\bs{x})}u_{ij}
\end{align}
That is,  $j_i^r$ is the available option with the highest realized reward and $j_i^u$ is the available option with the highest realized utility. Then, we define set $\Gamma(\bs{x})$ as the set of \textit{cooperative} scenarios where customer choices align with the planner's interests, i.e.,
\begin{align}
\Gamma(\bs{x}) = \{i\in\mathcal{I}:j_i^r = j_i^u\}
\end{align}
Then, Constraint (\ref{constr:SHP}) reduces to
\begin{equation}
\sum_{j\in \ml{J}:u_{ij}<u_{ik}}y_{ij} + x_k \le 1, \forall i \in \mathcal{I}\setminus \Gamma(\bs{x}),\forall k \in \mathcal{J}. \label{constr:rSHP}
\end{equation}
Clearly,  for MSMFLP, we have $\Gamma(\bs{x})= \mathcal{I},\forall \bs{x}\in\Omega$, and Constraint (\ref{constr:SHP}) is completely redundant. Consequently, MSMFLP reduces to the following model
\begin{subequations}\label{model:SAFLP}
\begin{alignat}{2}
\label{SABF:obj}  \max~~ & \frac{1}{N}\sum_{i\in \mathcal{I}}\sum_{j\in \mathcal{J}}r_{ij}y_{ij}\\
\label{SABF:cons1}\text{s.t.}~~& (\ref{SHP:constr0}), (\ref{SHP:constr1}), (\ref{SHP:constr2}), (\ref{SHP:constr4})
\end{alignat}
\end{subequations}
whose structure resembles that of the \textit{uncapacitated facility location problem} or the \textit{p-median problem}, neither of which involves customer choice-related constraints in their formulation. 
To the best of our knowledge, this paper is the first to discover the redundancy of choice-related constraints in choice-based facility location problems employing the partially binary rule.

Unlike MSMFLP, the other two application problems do not exhibit the complete redundancy of  Constraint (\ref{constr:SHP}). To further investigate how the reward-utility relationship influences the redundancy level of   Constraint (\ref{constr:SHP}) and the overall computational difficulty of \textbf{[ADPP]}, we conduct additional experiments on CAOP under the exponomial choice model, utilizing instances with varying reward-utility correlations.

To this end, we first generate the reward $R_j, \forall j\in\ml{J},$ using a lognormal distribution with a location parameter 0 and a scale parameter 0.2.  We then assume that the utility of an option $j \in \ml{J}$ is $U_j = V_j - \xi_j$, with the random component $\xi_j$ following  an exponential distribution $\xi_j \sim exp(1)$. The deterministic  component of the utility $V_j$ is defined as $V_j = \kappa R_j + b_j,\forall j\in\mathcal{J}$, where $b_j$ is a fixed parameters drawn  from the normal distribution with mean 1 and variance 1. The coefficient $\kappa$ quantifies  to what extent $V_j$ is related to $R_j$. For the outside option, we set $R_0 = V_0 = 0$.

\begin{table}[htb]
\centering
\caption{Impact of $\kappa$ on computational difficulty. Each row shows the summary statistics over 10 instances.}
\label{tab:difficulty}
\resizebox{0.8\textwidth}{!}{%
\begin{tabular}{lcccccccccccccccc}
\toprule[1.5pt]
$\kappa$ &  & \multicolumn{3}{c}{t{[}s{]}} &  & \multicolumn{3}{c}{\#Node} &  & \multicolumn{3}{c}{rgap{[}\%{]}} &  & \multicolumn{3}{c}{$\Gamma(\bs{x^*}) / N$} \\ \cline{3-5} \cline{7-9} \cline{11-13} \cline{15-17}
   &  & min  & max    & avg   &  & min & max  & avg     &  & min  & max  & avg  &  & min  & max  & avg  \\ \hline
-1 &  & 5.91 & 121.41 & 46.50 &  & 693 & 4480 & 2475.70 &  & 1.39 & 8.21 & 3.17 &  & 0.00 & 0.82 & 0.29 \\
0  &  & 0.91 & 52.37  & 14.98 &  & 9   & 2447 & 875.50  &  & 0.06 & 2.73 & 1.08 &  & 0.00 & 0.92 & 0.51 \\
1  &  & 0.88 & 13.37  & 6.21  &  & 11  & 739  & 352.30  &  & 0.01 & 1.41 & 0.70 &  & 0.04 & 0.99 & 0.69 \\ \bottomrule[1.5pt]
\end{tabular}%
}
\end{table}

We adjust $\kappa$ within $\{-1,0,1\}$ to depict three relations: \textit{negatively related}, \textit{not related} and \textit{positively related}. We maintain a constant setup with $|\mathcal{J}| = 100$ options and $N = 300$ samples (scenarios). For each $\kappa$ value, we change random seeds to generate 10 instances. We summarize the average values in Table \ref{tab:difficulty}, where $\Gamma(\bs{x^*}) / N$ is the proportion of scenarios rendered redundant of Constraint (\ref{constr:SHP}) by the optimal solution $\bs{x^*}$. Obviously, with the increase of $\kappa$, the value of $\Gamma(\bs{x^*}) / N$ increases. Meanwhile, we observe a decline in the root gap, the number of explored nodes, and the computational time. These findings suggest that the computational difficulty  of \textbf{[ADPP]} diminishes with the increase of $\kappa$. That is, \textbf{[ADPP]} is easier to solve when the reward and utilities exhibit a more positive correlation.

\eat{
Figure \ref{fig:difficulty} gives a more straightforward comparison of computational time. With the increase of $k$, we clearly see the reduction of computational time as well as time variance.
\begin{figure}[htb]
		\begin{center}
			\subfigure[Number of instances solved versus computational time]{
				\psfig{figure = figures/complexity, width=77mm, height = 55mm}\label{fig:app2-cum}}
			\subfigure[Boxplot of computational time]{
				\psfig{figure = figures/complexity-bp, width=77mm, height = 55mm}\label{fig:app2-bp} }
		\end{center}
		\caption{Time comparison of different $k$ values.}
		\label{fig:difficulty}
\end{figure}
}

\section{Concluding remarks}\label{sec:conclusion}

This paper studied a general class of choice-based DPP under uncertainties. We formulated a unified and tight MILP approximation for the general problem using sample average approximation. Three application problems were given to showcase the applicability of our framework. However, the formulation increases problem size, especially when the sample size is large. To efficiently solve our model, we proposed a Benders decomposition framework with customized cut separation procedures. Through extensive computational experiments, we demonstrated the solution quality and efficiency of our framework. Specifically, we showed that 1) our framework can obtain high-quality solutions with a reasonable number of scenarios; and 2) our framework outperforms advanced solution approaches in the literature among three applications. Finally, we use CAOP as an example to reveal the impact of reward-utility relation on computational difficulty and derived valuable computational insights.

Our Benders decomposition approach is well-suited for cases where the uncertainty $\bs{\xi}$ follows a continuous distribution. Under such a condition, the choice variable $\bs{y}$ can be relaxed from binary to continuous,  facilitating the effective development of our algorithm. However, when $\bs{\xi}$ follows a discrete distribution or when limited distribution information is available, the choice variable $\bs{y}$ may not be relaxed without the loss of optimality. Consequently, the conventional Benders method is not readily applicable, underscoring the necessity for algorithm adjustments and improvements. \JJ{A notable extension for choice-based DPP involves integrating resource capacity constraints for each option (be it a service or facility) by accounting for the stochastic nature of customer arrivals, such as through a Poisson process. Given that the planner does not control customer choice behaviors, employing a chance-constrained program (CCP) could effectively model capacity constraints. Nonetheless, this adaptation makes the subproblem non-separable and the choice variable $\bs{y}$ cannot be relaxed, potentially limiting the use of our Benders decomposition approach. We would like to leave the model and algorithm development for future research.}

\JJ{To handle uncertainty, our framework falls into the stochastic program regime. We show that our framework works very well under a reasonable number of samples. However, there are also cases where only very limited samples are available. As such, robust optimization techniques may also bring new modeling and managerial insights for discrete planning problems under uncertainty. For example, given a set of samples as nominal distribution, a distributional robust model based on Wasserstein distance can be built to mitigate the risk caused by limited sample size. Efforts need to be paid to find tractable robust counterpart and efficient solution approaches.}
\theendnotes

\section*{Acknowledgement}
The authors would like to thank Ali Aouad for his generosity in sharing the data and the code for the assortment optimization problem under the exponomial choice model. 

\bibliographystyle{apalike}

 \let\oldbibliography\thebibliography
 \renewcommand{\thebibliography}[1]{%
  \oldbibliography{#1}%
  \baselineskip 11pt 
  \setlength{\itemsep}{5pt}
 }

\bibliography{ref}

\vspace{15mm}

\begin{center}
\Large Appendix for ``Approximate Resolution of Stochastic Choice-based Discrete Planning''
\end{center}

\section*{A.1.~Proofs}
\subsection*{A.1.1.~Proof of Proposition 1}

Given a feasible $\bs{\bar{x}}\in\Omega$, define $\mathcal{J}^o = \left\{j\in\mathcal{J}: \bar{x}_j = 1\right\}$ and $j_i^* = \arg\max_{j'\in \mathcal{J}^o}u_{ij'}$.  An optimal primal solution $\bar{y}$ to \textbf{[SP$_i$]} can be easily obtained as
\begin{align}\tag{A.1}
\bar{y}_{ij}=\left\{\begin{matrix}
1,  & j = j_i^*\\
0,  & \text{otherwise}
\end{matrix}\right.
\end{align}
Now, given $\bs{\bar{x}}$ and $\bs{\bar{y}}$, the reduced KKT conditions without primal feasibility  are
\begin{align}
&\bar{\nu}_{ij}(\bar{x}_j-\bar{y}_{ij})=0,\forall i\in \mathcal{I},\forall j\in \mathcal{J},\label{cs1}\tag{A.2}\\
&\bar{\mu}_{ik}(1 - \bar{x}_k - \sum_{j\in \ml{J}:u_{ij}<u_{ik}}\bar{y}_{ij})=0,\forall i\in \mathcal{I},\forall k\in \mathcal{J},\label{cs2}\tag{A.3}\\
&r_{ij} - \bar{\lambda}_i - \bar{\nu}_{ij} - \sum_{k\in \ml{J}:u_{ik}>u_{ij}}\bar{\mu}_{ik} \le 0,\forall i\in \mathcal{I},\forall j\in \mathcal{J},\label{sc}\tag{A.4}\\
&\bar{\bs{\nu}}_i,\bar{\bs{\mu}}_i \ge 0.\tag{A.5}
\end{align}

When $j\in\mathcal{J}^o$, we have $\bar{\nu}_{ij}(\bar{x}_j-\bar{y}_{ij})=0$ since $\bar{\nu}_{ij} = 0$; when $j\in\mathcal{J}\setminus \mathcal{J}^o$, we have $\bar{\nu}_{ij}(\bar{x}_j-\bar{y}_{ij})=0$ since $\bar{x}_j - \bar{y}_{ij} = 0$. Therefore, (\ref{cs1}) is satisfied.

When $k = j_i^*$, we have $\bar{\mu}_{ik}(1 - \bar{x}_k - \sum_{j\in \ml{J}:u_{ij} < u_{ik}}\bar{y}_{ij})=0$ since $1 - \bar{x}_{j_i^*} - \sum_{j\in\mathcal{J}:u_{ij}<u_{ik}}\bar{y}_{ij} = 1 - 1 - 0 = 0$; when $k\in \mathcal{J}\setminus \left\{j_i^*\right\}$, we have $\bar{\mu}_{ik}(1 - \bar{x}_k - \sum_{j\in \ml{J}:u_{ij} < u_{ik}}\bar{y}_{ij})=0$ since $\bar{\mu}_{ik} = 0$. Therefore, (\ref{cs2}) is satisfied.

When $j = j_i^*$, we have $r_{ij} - \bar{\lambda}_i - \bar{\nu}_{ij} - \sum_{k\in \ml{J}:u_{ik} > u_{ij}}\bar{\mu}_{ik} = r_{ij_i^*} - r_{ij_i^*} - 0 - 0 = 0$; when $j\in\mathcal{J}^o\setminus\left\{j_i^*\right\}$, we have $r_{ij} - \bar{\lambda}_i - \bar{\nu}_{ij} - \sum_{k\in \ml{J}:u_{ik}>u_{ij}}\bar{\mu}_{ik} = r_{ij} - \bar{\lambda}_i - [\max_{j\in \mathcal{J}^o\setminus \left\{j_i^*\right\}   }r_{ij}-\bar{\lambda}_{i}]_{+}\le 0$; when $j\in \mathcal{J}\setminus \mathcal{J}^o$, we have $r_{ij} - \bar{\lambda}_i - \bar{\nu}_{ij} - \sum_{k\in \ml{J}:u_{ik} > u_{ij}}\bar{\mu}_{ik} = r_{ij} - \bar{\lambda}_i - \sum_{k\in \ml{J}:u_{ik}>u_{ij}}\bar{\mu}_{ik} - [r_{ij}-\bar{\lambda}_{i} - \sum_{k\in \mathcal{J}:u_{ik}>u_{ij}}\bar{\mu}_{ik}]_{+} \le 0$. Therefore, (\ref{sc}) is satisfied.

Finally, we check the dual and primal objectives, i.e., $\bar{\lambda}_i + \sum_{j \in \mathcal{J}} \bar{\nu}_{ij} \bar{x}_{j} + \sum_{k\in \mathcal{J}}\bar{\mu}_{ik}(1-\bar{x}_k )
=  r_{ij_i^*} + \sum_{j\in\mathcal{J}^o}0\cdot 1 + \sum_{j\in\mathcal{J}\setminus\mathcal{J}^o}[r_{ij}-\bar{\lambda}_{i} - \sum_{k\in \mathcal{J}:u_{ik}>u_{ij}}\bar{\mu}_{ik}]_{+} \cdot 0 + 0\cdot (1-1) + \sum_{j\in\mathcal{J}\setminus \left\{j_i^*\right\}}0\cdot (1-0) = r_{ij_i^*} = \sum_{j\in\mathcal{J}}r_{ij}\bar{y}_{ij}$. That is, the primal objective is equal to the dual objective. Altogether, we show that the proposed ($\bar{\lambda}_i, \bar{\bs{\nu}}_i, \bar{\bs{\mu}}_i$) is an optimal solution to \textbf{[DSP$_i$]}.\quad$\square$

\subsection*{A.1.2.~Proof of Proposition 2}

To start, Constraint (25c) can be written as
\begin{align}
&y_{ij} \le x_j,\forall j\in\mathcal{J},\tag{A.6}\\
&y_{ij} \le 1 - \delta_{ijk}x_k,\forall j\in\mathcal{J},\forall k\in\mathcal{J} \label{constr:nCAC}.\tag{A.7}
\end{align}
It suffices to show that Constraint (\ref{constr:nCAC}) has the same effect as Constraint (15e). For any $i\in\mathcal{I}, k\in\mathcal{J}$, define set $\mathcal{J}_i^k = \left\{j\in \mathcal{J}:U_{ij} < U_{ik}\right\}$. Then Constraint (\ref{constr:nCAC}) is equivalent to $y_{ij} \le 1 - x_k,\forall j\in \mathcal{J}_i^k$. If option $k$ is offered, then we have $x_k = 1$ and $y_{ij} = 0,\forall j\in \mathcal{J}_i^k$, meaning that any option with lower utility than option $k$ will not be chosen; if option $k$ is not offered, then we have $x_k = 0$ and $y_{ij}\le 1,\forall j\in\mathcal{J}_i^k$, which is automatically satisfied. Therefore, Constraint (\ref{constr:nCAC}) has the same effect as Constraint (15e) and \textbf{[nSP$_i$]} is equivalent to \textbf{[SP$_i$]}.\quad $\square$

\subsection*{A.1.3.~Proof of Proposition 3}

Under each scenario $i\in\mathcal{I}$, define $\mathcal{J}^c = \left\{j\in \mathcal{J}: r_{ij} > r_{ij_i^*}\right\}$. Then we can construct $\bs{\bar{y}}$ as
\begin{align}\tag{A.8}
\bar{y}_{ij}=\left\{\begin{matrix}
\bar{\beta}_{ij}  &,\forall j\in \mathcal{J}^c, \\
1-\sum_{j\in \mathcal{J}^c}\bar{\beta}_{ij}  &,j=j_i^*, \\
0  &,\text{otherwise}
\end{matrix}\right.
\end{align}
The reduced KKT conditions without primal feasibility  are
\begin{align}
&r_{ij} - \bar{\lambda}'_i-\bar{\eta}_{ij} \le 0,\forall j\in J,\label{KKT:sc}\tag{A.9}\\
&\bar{\eta}_{ij}(\bar{\beta}_{ij} - \bar{y}_{ij}) = 0,\forall j\in J,\label{KKT:cs}\tag{A.10}\\
&\bar{\eta}_{ij}\ge 0,\forall j\in J.\tag{A.11}
\end{align}

We thus have $\bar{\lambda}'_i + \bar{\eta}_{ij}  = r_{ij_i^*} +[r_{ij} - r_{ij_i^*}]_{+}  \ge r_{ij},\forall j\in\mathcal{J}$. Therefore, (\ref{KKT:sc}) is satisfied.

When $j\in\mathcal{J}^c$, we have $\bar{\eta}_{ij}(\bar{\beta}_{ij} - \bar{y}_{ij}) = 0$ since $\bar{\beta}_{ij} - \bar{y}_{ij} = 0$; when $j = \mathcal{J}\setminus \mathcal{J}^c$, we have $\bar{\eta}_{ij}(\bar{\beta}_{ij} - \bar{y}_{ij}) = 0$ since $\bar{\eta}_{ij} = 0$. Therefore, (\ref{KKT:cs}) is satisfied.

Finally, we check the primal and dual objectives:  $\bar{\lambda}'_i + \sum_{j\in \mathcal{J}}\bar{\eta}_{ij}\bar{\beta}_{ij}  = \bar{\lambda}'_i +\sum_{j\in \mathcal{J}^c}\bar{\eta}_{ij}\bar{y}_{ij} + \sum_{j\in \mathcal{J}\setminus \mathcal{J}^c}\bar{\eta}_{ij}\bar{\beta}_{ij}
=  r_{ij_i^*} + \sum_{j\in \mathcal{J}^c}(r_{ij} - r_{ij_i^*})\bar{y}_{ij} =  r_{ij_i^*} (1-\sum_{j\in \mathcal{J}^c}\bar{y}_{ij}) + \sum_{j\in \mathcal{J}^c}r_{ij}\bar{y}_{ij} = \sum_{j\in \mathcal{J}}r_{ij}\bar{y}_{ij}$,
which states that dual and primal objectives are equal. Altogether, we prove that the proposed $(\bar{\lambda}'_i,\bs{\bar{\eta}_i})$ is an optimal solution to \textbf{[nDSP$_i$]}.\quad $\square$

\subsection*{A.1.4.~Proof of Proposition 4}

Let $k_{ij}^* = \arg\min_{k\in \mathcal{J}} 1-\delta_{ijk}\bar{x}_k$. It directly follows from the definition of $\bs{\beta}, \bs{\bar{\nu}}', \bs{\bar{\mu}}'$ that
\begin{align}
\theta_i \le &\bar{\lambda}'_i + \sum_{j\in \mathcal{J}}\bar{\eta}_{ij}\bar{\beta}_{ij}\tag{A.13}\\
= &\bar{\lambda}'_i + \sum_{j\in \mathcal{J}}\bar{\eta}_{ij}\left(\mathbb{I}\left \{ \bar{\beta}_{ij} = \bar{x}_j \right \}\cdot x_{j} + \mathbb{I}\left \{ \bar{\beta}_{ij} = 1 - \delta_{ijk^*_{ij}}\bar{x}_{k^*_{ij}} \right \}\cdot (1 - \delta_{ijk^*_{ij}}x_{k^*_{ij}})\right)\tag{A.14}\\
= & \bar{\lambda}_i' + \sum_{j\in \mathcal{J}}\bar{\nu}_{ij}'x_j  + \sum_{j\in \mathcal{J}}\bar{\mu}'_{ijk_{ij}^*} (1- \delta_{ijk_{ij}^*}x_{k_{ij}^*}) \label{cut:frac1}\tag{A.15}\\
= & \bar{\lambda}_i' + \sum_{j\in \mathcal{J}}\bar{\nu}_{ij}'x_j  + \sum_{j\in \mathcal{J}}\sum_{k\in\mathcal{J}}\bar{\mu}'_{ijk} (1- \delta_{ijk}x_k) \tag{A.16}\label{cut:frac2}\\
= &  \bar{\lambda}_i' + \sum_{j\in \mathcal{J}}\sum_{k\in\mathcal{J}}\bar{\mu}'_{ijk} + \sum_{j\in \mathcal{J}}(\bar{\nu}'_{ij} - \sum_{k\in J}\delta_{ijk}\bar{\mu}'_{ijk})x_j
\tag{A.17}
\end{align}
 Equalities (\ref{cut:frac1}) and (\ref{cut:frac2}) comes from the fact that $\bar{\mu}'_{ijk} = 0,\forall k\ne k_{ij}^*,\forall i\in\mathcal{I},\forall j\in\mathcal{J}$. \quad $\square$

\section*{A.2.~Solution quality test}

This section provides solution quality results on three applications. We also compare Monte Carlo Sampling (MCS) and Latin Hypercube Sampling (LHS). 

\subsection*{A.2.1.~Results of  CAOP under multinomial probit model}
We apply our SBBD on the multinomial probit model, which does not have closed-form formulation of the choice probability. Here, the utility of an alternative $j\in \ml{J}$  is given by $U_j = V_j + \xi_j$, where $V_j$ a deterministic component and $\xi_j$ follow some i.i.d. normal distribution. We generate the problem instances as follow: both the deterministic utility $V_j$ and reward $R_j$ is generated from a uniform distribution on the interval [0,100] (with $R_0 =0$). The random utility $\xi$ follows normal distribution with mean 0 and variance 100. 

We use the estimated gaps $\Delta$ and $\Delta_{0.95}$ to assess the solution quality. Table \ref{tab:app1-sol} reports the results. Obviously, LHS reduces both the variance $\sigma$ and gaps compared with MCS. Moreover, when $N=1000$, our SBBD supported by LHS is able to obtain small $\Delta$ values and thus proven high-quality solution.
\begin{table}[h]
\small
\centering
\caption{Solution validation for CAOP under  multinomial probit model.}
\label{tab:app1-sol}
\resizebox{0.8\textwidth}{!}{%
\begin{tabular}{ccccccccccccccc}
\toprule[1.5pt]
\multirow{2}{*}{$N$} &
  \multirow{2}{*}{$|\mathcal{A}|$} &
  \multirow{2}{*}{$\tau$} &
   &
  \multicolumn{5}{c}{MCS} &
   &
  \multicolumn{5}{c}{LHS} \\ \cline{5-9} \cline{11-15} 
 &
   &
   &
   &
  $\hat{v}$ &
  $\bar{v}^M_N$ &
  $\sigma$ &
  $\Delta$ &
  $\Delta_{0.95}$ &
   &
  $\hat{v}$ &
  $\bar{v}^M_N$ &
  $\sigma$ &
  $\Delta$ &
  $\Delta_{0.95}$ \\ \hline
300  & 100 & 10 &  & 82.48 & 83.88 & 0.3853 & 1.70 & 2.47 &  & 82.45 & 83.39 & 0.2863 & 1.14 & 1.71 \\
500  & 100 & 10 &  & 82.77 & 83.83 & 0.4011 & 1.28 & 2.08 &  & 82.74 & 83.16 & 0.2626 & 0.51 & 1.03 \\
1000 & 100 & 10 &  & 82.95 & 83.26 & 0.4966 & 0.37 & 1.36 &  & 82.83 & 83.05 & 0.0943 & 0.27 & 0.45 \\
300  & 200 & 10 &  & 86.57 & 88.12 & 0.4224 & 1.79 & 2.60 &  & 86.83 & 88.10 & 0.4275 & 1.46 & 2.27 \\
500  & 200 & 10 &  & 87.03 & 88.35 & 0.4407 & 1.52 & 2.35 &  & 87.14 & 87.97 & 0.4138 & 0.95 & 1.74 \\
1000 & 200 & 10 &  & 87.45 & 88.33 & 0.5272 & 1.01 & 2.00 &  & 87.52 & 87.89 & 0.3809 & 0.42 & 1.14 \\ \bottomrule[1.5pt]
\end{tabular}%
}
\end{table}

\subsection*{A.2.2.~Results of FLoP and MSMFLP}

This section provides complete results of solution quality test for FLoP and MSMFLP. The dataset is the same as Section 5 in the main text. Tables \ref{tab:app2-sol} and \ref{tab:app3-sol} report the estimated gaps.  Obviously, LHS significantly reduces the variance of objective as well as the approximation error compared with MCS. Moreover, for LHS, we include the change of $\Delta$ and $\Delta_{0.95}$ when problem size (i.e., $|\ml{A}|$ and $|\ml{J}|$) increases in the parenthesis. For example, given $N = 300$ and $\tau = 5$, $\Delta$ increases from 0.83\% to 1.79\% when $|\ml{A}|$ increases from 10 to 20. Therefore, we put $1.79-0.83 = 0.96$ in the parenthesis. Overall, we observe only slight increase of $\Delta$ and $\Delta_{0.95}$ for both FLoP (less than 1\%)  and MSMFLP (less than 0.4\%).

\begin{table}[htb]
\centering
\caption{Solution validation for FLoP.}
\label{tab:app2-sol}
\resizebox{0.8\textwidth}{!}{%
\begin{tabular}{ccccccccccccccc}
\toprule[1.5pt]
\multirow{2}{*}{$N$} &
  \multirow{2}{*}{$|\mathcal{A}|$} &
  \multirow{2}{*}{$\tau$} &
   &
  \multicolumn{5}{c}{MCS} &
   &
  \multicolumn{5}{c}{LHS} \\ \cline{5-9} \cline{11-15} 
 &
   &
   &
   &
  $\hat{v}$ &
  $\bar{v}^M_N$ &
  $\sigma$ &
  $\Delta$ &
  $\Delta_{0.95}$ &
   &
  $\hat{v}$ &
  $\bar{v}^M_N$ &
  $\sigma$ &
  $\Delta$ &
  $\Delta_{0.95}$ \\ \hline
300  & 10 & 5 &  & 3.79 & 3.90 & 1.7146 & 2.92 & 3.66 &  & 3.79 & 3.82 & 0.5632 & 0.83 & 1.08 \\
500  & 10 & 5 &  & 3.79 & 3.86 & 1.2640 & 1.94 & 2.49 &  & 3.79 & 3.81 & 0.4640 & 0.67 & 0.87 \\
1000 & 10 & 5 &  & 3.79 & 3.82 & 0.7711 & 0.77 & 1.11 &  & 3.79 & 3.80 & 0.3453 & 0.18 & 0.33 \\
300  & 20 & 5 &  & 3.91 & 4.09 & 1.4531 & 4.61 & 5.22 &  & 3.93 & 4.00 & 0.4413 & 1.79(+0.96) & 1.97(+0.89) \\
500  & 20 & 5 &  & 3.93 & 4.03 & 1.1287 & 2.48 & 2.96 &  & 3.93 & 3.98 & 0.4116 & 1.18(+0.51) & 1.35(+0.48) \\
1000 & 20 & 5 &  & 3.93 & 3.99 & 0.7212 & 1.37 & 1.68 &  & 3.93 & 3.95 & 0.3280 & 0.43(+0.25) & 0.57(+0.24) \\ \bottomrule[1.5pt]
\end{tabular}%
}
\end{table}

\begin{table}[htb]
\centering
\caption{Solution validation for MSMFLP.}
\label{tab:app3-sol}
\resizebox{0.8\textwidth}{!}{%
\begin{tabular}{ccccccccccccccc}
\toprule[1.5pt]
\multirow{2}{*}{$N$} &
  \multirow{2}{*}{$|\mathcal{J}|$} &
  \multirow{2}{*}{$\tau$} &
   &
  \multicolumn{5}{c}{MCS} &
   &
  \multicolumn{5}{c}{LHS} \\ \cline{5-9} \cline{11-15} 
 &
   &
   &
   &
  $\hat{v}$ &
  $\bar{v}^M_N$ &
  $\sigma$ &
  $\Delta$ &
  $\Delta_{0.95}$ &
   &
  $\hat{v}$ &
  $\bar{v}^M_N$ &
  $\sigma$ &
  $\Delta$ &
  $\Delta_{0.95}$ \\ \hline
300 &
  100 &
  20 &
   &
  38.18 &
  41.29 &
  0.2214 &
  8.15 &
  9.10 &
   &
  39.26 &
  39.86 &
  0.0854 &
  1.53 &
  1.89 \\
500 &
  100 &
  20 &
   &
  38.83 &
  39.42 &
  0.1636 &
  6.15 &
  6.85 &
   &
  39.42 &
  39.77 &
  0.0513 &
  0.89 &
  1.10 \\
1000 &
  100 &
  20 &
   &
  39.24 &
  40.35 &
  0.1252 &
  2.82 &
  3.35 &
   &
  39.46 &
  39.61 &
  0.0336 &
  0.38 &
  0.52 \\
300 &
  200 &
  20 &
   &
  38.80 &
  42.24 &
  0.2424 &
  8.87 &
  9.90 &
   &
  40.12 &
  40.85 &
  0.0892 &
  1.82(+0.29) &
  2.19(+0.30) \\
500 &
  200 &
  20 &
   &
  39.33 &
  41.89 &
  0.1857 &
  6.51 &
  7.29 &
   &
  40.19 &
  40.63 &
  0.0566 &
  1.10(+0.21) &
  1.33(+0.23) \\
1000 &
  200 &
  20 &
   &
  39.76 &
  41.69 &
  0.1184 &
  4.84 &
  5.53 &
   &
  40.23 &
  40.53 &
  0.0327 &
  0.74(+0.36) &
  0.87(+0.35) \\ \bottomrule[1.5pt]
\end{tabular}%
}
\end{table}

\section*{A.3.~Efficiency test on CAOP under multinomial probit model}

This section provides additional efficiency test results on CAOP{ under multinomial probit model.
The utility of an alternative $a\in \ml{A}$ is $U_a = V_a + \xi_a$.  Both $V_a$ and $R_a$ are generated from $\mathcal{U}[0,100], \forall a\in \mathcal{A}$.  We then set $V_0 = 50$ and $R_0 = 0$. The random term $\xi_j$ follows i.i.d. normal distribution with mean 0 and variance $\mathcal{V}$. We employ the same stabilization procedure in Stage 1 as the exponomial choice model in Section 5.1 of the main text.

Table \ref{tab:efficiency-probit} reports the computational result. Apparently, SBBD outperforms BIBC by a large margin. Moreover, compared to MILP, SBBD is faster for all the 24 instances. Particularly, the superiority of SBBD becomes even more significant for the eight  ``hard'' instances with $N = 1000$.

\begin{table}[htb]
\centering
\caption{Efficiency test on CAOP under multinomial probit model}\label{tab:efficiency-probit}
\resizebox{\textwidth}{!}{%
\begin{tabular}{cccccccccclcccclccccc}
\toprule[1.5pt]
\multirow{2}{*}{$N$} &
  \multirow{2}{*}{$|\mathcal{J}|$} &
  \multirow{2}{*}{$\tau$} &
  \multirow{2}{*}{$\mathcal{V}$} &
   &
  \multicolumn{5}{c}{SBBD} &
   &
  \multicolumn{4}{c}{MILP} &
   &
  \multicolumn{5}{c}{BIBC} \\ \cline{6-10} \cline{12-15} \cline{17-21} 
 &
   &
   &
   &
   &
  t{[}s{]} &
  \#N &
  \#C &
  rgap{[}\%{]} &
  ogap{[}\%{]} &
   &
  t{[}s{]} &
  \#N &
  rgap{[}\%{]} &
  ogap{[}\%{]} &
   &
  t{[}s{]} &
  \#N &
  \#C &
  rgap{[}\%{]} &
  ogap{[}\%{]} \\ \hline
300  & 100 & 10 & 100 &  & 6.2  & 85   & 10666 & 0.28 & 0.00 &  & 7.4   & 1  & 0.00 & 0.00 &  & 142.4  & 3984  & 2513 & 3.80 & 0.00 \\
     &     &    & 200 &  & 4.7  & 38   & 8015  & 0.38 & 0.00 &  & 10.0  & 1  & 0.00 & 0.00 &  & 145.1  & 4329  & 2184 & 3.49 & 0.00 \\
     &     & 20 & 100 &  & 5.0  & 1    & 9092  & 0.00 & 0.00 &  & 6.7   & 1  & 0.00 & 0.00 &  & 727.3  & 30937 & 2266 & 5.47 & 0.00 \\
     &     &    & 200 &  & 5.0  & 11   & 8929  & 0.07 & 0.00 &  & 6.9   & 1  & 0.00 & 0.00 &  & 2411.0 & 57162 & 2159 & 5.14 & 0.00 \\
     & 200 & 10 & 100 &  & 14.4 & 1266 & 15554 & 1.39 & 0.00 &  & 43.1  & 24 & 0.92 & 0.00 &  & 498.4  & 10274 & 3230 & 2.92 & 0.00 \\
     &     &    & 200 &  & 10.9 & 550  & 11090 & 1.06 & 0.00 &  & 44.4  & 49 & 0.73 & 0.00 &  & 214.3  & 4921  & 3126 & 2.01 & 0.00 \\
     &     & 20 & 100 &  & 10.7 & 304  & 10885 & 0.63 & 0.00 &  & 63.8  & 9  & 0.36 & 0.00 &  & 3600.0 & 39168 & 3576 & 3.73 & 1.45 \\
     &     &    & 200 &  & 9.8  & 171  & 11357 & 0.50 & 0.00 &  & 64.9  & 7  & 0.22 & 0.00 &  & 3600.0 & 31076 & 3694 & 3.01 & 1.14 \\
500  & 100 & 10 & 100 &  & 10.2 & 1    & 23464 & 0.00 & 0.00 &  & 10.0  & 1  & 0.00 & 0.00 &  & 443.6  & 5196  & 4564 & 4.41 & 0.00 \\
     &     &    & 200 &  & 8.3  & 203  & 12840 & 0.70 & 0.00 &  & 11.9  & 1  & 0.00 & 0.00 &  & 195.0  & 4322  & 3602 & 2.94 & 0.00 \\
     &     & 20 & 100 &  & 4.7  & 1    & 7630  & 0.00 & 0.00 &  & 10.0  & 1  & 0.00 & 0.00 &  & 3600.0 & 32683 & 3860 & 5.79 & 0.52 \\
     &     &    & 200 &  & 7.4  & 19   & 10769 & 0.25 & 0.00 &  & 10.1  & 1  & 0.00 & 0.00 &  & 3600.0 & 23834 & 3763 & 5.33 & 0.98 \\
     & 200 & 10 & 100 &  & 21.6 & 732  & 22384 & 1.16 & 0.00 &  & 106.7 & 36 & 0.59 & 0.00 &  & 1338.5 & 7818  & 4025 & 3.20 & 0.00 \\
     &     &    & 200 &  & 27.6 & 1125 & 24099 & 1.27 & 0.00 &  & 141.9 & 29 & 0.59 & 0.00 &  & 3600.0 & 18919 & 4751 & 3.95 & 0.87 \\
     &     & 20 & 100 &  & 23.0 & 538  & 22674 & 0.54 & 0.00 &  & 131.4 & 1  & 0.00 & 0.00 &  & 3600.0 & 10352 & 5187 & 2.39 & 0.69 \\
     &     &    & 200 &  & 22.5 & 446  & 22393 & 0.41 & 0.00 &  & 107.2 & 15 & 0.55 & 0.00 &  & 3600.0 & 9691  & 4999 & 3.08 & 1.58 \\
1000 & 100 & 10 & 100 &  & 22.9 & 124  & 41533 & 0.38 & 0.00 &  & 49.1  & 1  & 0.00 & 0.00 &  & 1554.1 & 2742  & 7234 & 3.85 & 0.00 \\
     &     &    & 200 &  & 29.4 & 501  & 41781 & 0.81 & 0.00 &  & 44.0  & 1  & 0.00 & 0.00 &  & 2001.9 & 4447  & 9963 & 3.76 & 0.00 \\
     &     & 20 & 100 &  & 17.4 & 1    & 26118 & 0.00 & 0.00 &  & 51.4  & 1  & 0.00 & 0.00 &  & 3600.0 & 5284  & 7544 & 5.39 & 2.46 \\
     &     &    & 200 &  & 15.2 & 1    & 23568 & 0.00 & 0.00 &  & 52.0  & 1  & 0.00 & 0.00 &  & 3600.0 & 3726  & 6800 & 5.36 & 2.68 \\
     & 200 & 10 & 100 &  & 59.5 & 1285 & 33294 & 1.91 & 0.00 &  & 228.3 & 15 & 0.89 & 0.00 &  & 3600.0 & 2982  & 8859 & 2.44 & 0.87 \\
     &     &    & 200 &  & 83.7 & 1515 & 31773 & 1.95 & 0.00 &  & 341.8 & 67 & 1.21 & 0.00 &  & 3600.0 & 2855  & 8539 & 2.76 & 1.44 \\
     &     & 20 & 100 &  & 52.4 & 687  & 50238 & 0.52 & 0.00 &  & 293.4 & 18 & 0.77 & 0.00 &  & 3600.0 & 2225  & 9867 & 3.10 & 2.62 \\
     &     &    & 200 &  & 60.3 & 537  & 58418 & 0.59 & 0.00 &  & 336.5 & 47 & 0.98 & 0.00 &  & 3600.0 & 1649  & 9902 & 2.93 & 2.43 \\ \bottomrule[1.5pt]
\end{tabular}%
}
\end{table}

\end{document}